\documentclass[journal,onecolumn,draftclsnofoot,]{IEEEtran}

\usepackage{setspace}

\usepackage{enumerate}

\usepackage[]{graphicx,float,latexsym,times}
\usepackage{amsfonts,amstext,amsmath,amssymb,amsthm}
\usepackage{caption2}
\usepackage{multirow}
\usepackage{subfigure}
\usepackage{xcolor}
\usepackage{algorithm}
\usepackage{mathtools}
\usepackage{enumitem}
\usepackage{amsfonts,amssymb}
\usepackage{xcolor}
\usepackage{hyperref}
\usepackage[noend]{algpseudocode}
\usepackage{pdflscape}
\usepackage{booktabs}
\usepackage{enumerate}
\usepackage{enumitem}
\usepackage{stackrel}
\usepackage{lipsum}

\makeatletter
\newcommand{\distas}[1]{\mathbin{\overset{#1}{\kern\z@\sim}}}%
\newsavebox{\mybox}\newsavebox{\mysim}
\newcommand{\distras}[1]{%
	\savebox{\mybox}{\hbox{\kern3pt$\scriptstyle#1$\kern3pt}}%
	\savebox{\mysim}{\hbox{$\sim$}}%
	\mathbin{\overset{#1}{\kern\z@\resizebox{\wd\mybox}{\ht\mysim}{$\sim$}}}%
}
\makeatother

  \newtheorem{assumption}{Assumption}[section]
  \newtheorem{lemma}{Lemma}
\newtheorem{theorem}{Theorem}

\long\def\symbolfootnote[#1]#2{\begingroup%
\def\thefootnote{\fnsymbol{footnote}}\footnote[#1]{#2}\endgroup}

\newcommand{\hatQ}{{\widehat Q}}
\newcommand{\hatU}{{\widehat U}}
\newcommand{\hatG}{{\widehat G}}
\newcommand{\hatLambda}{{\widehat \Lambda}}

\ifCLASSINFOpdf
\else
\fi

\begin{document}
%
\title{Spectral CUSUM for Online Network Structure Change Detection}
%
%
%

\author{Minghe~Zhang,
        Liyan~Xie,
        and~Yao~Xie
\thanks{Minghe Zhang (Email:
	mzhang388@gatech.edu), and Yao Xie (Email:
	yao.xie@isye.gatech.edu, corresponding author) are with H. Milton Stewart School of Industrial and Systems Engineering, Georgia Institute of Technology, Atlanta,
GA, 30332 USA.}
\thanks{Liyan Xie (Email:
	xieliyan@cuhk.edu.cn) is with School of Data Science, The Chinese University of Hong Kong, Shenzhen, China.}
}

\markboth{Submitted to IEEE Transactions on Information Theory}%
{Shell \MakeLowercase{\textit{et al.}}: Bare Demo of IEEEtran.cls for IEEE Journals}
%



\maketitle

\begin{abstract}
Detecting abrupt changes in the community structure of a network from noisy observations is a fundamental problem in statistics and machine learning. This paper presents an online change detection algorithm called Spectral-CUSUM to detect unknown network structure changes through a generalized likelihood ratio statistic. We characterize the average run length (ARL) and the expected detection delay (EDD) of the Spectral-CUSUM procedure and prove its asymptotic optimality. Finally, we demonstrate the good performance of the Spectral-CUSUM procedure and compare it with several baseline methods using simulations and real data examples on seismic event detection using sensor network data.
\end{abstract}


%
\IEEEpeerreviewmaketitle

\section{Introduction}
\label{sec:intro}


Detecting network structure change from sequential data is a fundamental problem in high-dimensional data analysis, emerging from multiple applications, including seismic sensor networks \cite{faulkner2011next}, traffic networks \cite{ahmed2008novel}, swarm behavior monitoring \cite{berger2016classifying}, and social network change detection \cite{basuchowdhuri2019fast}. A community corresponds to a subset of nodes with much higher connectivity within the group than across groups. Real-world communication structure changes can be complicated. In various settings, the change may correspond to the emergence of a community, switching community memberships, changes in the number of communities, etc.

The need for community structure change detection is high for high-dimensional sequential data, which tend to have complex inter-dependent relationships between different dimensions. Such interdependence structure can be explicit, where the network topology be inferred from data. In the network settings, {\it the characteristic of the changes} will be a shift in {\it structures} of the underlying parameters, which is fundamentally different from a simple mean-shift considered in the change-point detection literature. As has been recognized, it is possible to exploit the underlying dependence structures to design asymptotically optimal detection algorithms \cite{keshavarz2018optimal}. In addition, we typically need a detection procedure to be computationally and memory efficient, as the data streams in this setting are often very high-dimensional and generated at high speed.

We first give a few examples of detecting changes exploring community structures:
\begin{itemize}
	\item Seismic sensor networks: The seismic research has been based on the massive amount of continuous data recorded by ultra-dense seismic sensor arrays, and many such data are publicly available on IRIS (\url{https://www.iris.edu}). In the old days, network seismology treated seismic signals individually - one sensor at a time - and detected an earthquake when multiple impulsive arrivals were consistent with a source within the Earth \cite{johnson1997robust}. Recently, with advances in sensor technology, which bring densely sampled data and high-performance computing and communication, we may be able to use a {\it network-based detection} by exploiting correlations between sensors to extract coherence signals. This will enhance the systematic detection of weak and unusual events that currently go undetected using individual sensors. Detecting such weak events is very crucial for earthquake prediction \cite{vere1970stochastic,omi2013forecasting}, oil field exploration, volcano monitoring, and deeper earth studies \cite{huang2015yellowstone}. 
	
	\item Social networks: The widespread use of social networks leads to a large amount of user-generated continuous data, which is quite valuable in studying many social phenomena. One important application is to detect change points using social network data. These change points may represent the collective anticipation of or response to external events or system ``shocks'' \cite{peel2015detecting}. Detecting such changes can provide a better understanding of patterns of social life. In other cases, early detection of change-point can predict or even prevent social stress due to disease or international threats. In social network data
	\cite{du2007community}, each node represents one individual, and the edge represents the relationship between two individuals. 
	
	\item Manifold is a common low-dimensional structure that lies in high-dimensional data, which can be captured using a similar network such as Isomap and Laplacian Eigenmaps \cite{talwalkar2008large}. Thus, change in the manifold structure can be detected from similar graphs. 
\end{itemize} 


In this paper, we present a new online change-point detection procedure, the Spectral CUSUM, to detect network structure changes by observing node features. We model the network structure through the inverse covariance of the noisy Gaussian features. Based on such a model, Spectral CUSUM is derived based on generalized likelihood ratios, where the unknown post-change parameters are estimated sequentially. This approach enables us to detect general types of structural changes, including the emergence of community and switching membership. The main theoretical contribution is to show the first-order asymptotic optimality of Spectral CUSUM and characterize the optimal choice of the parameters. We also present an online scheme for computing the detection statistic based on subspace tracking that is computationally and memory efficient. We demonstrate the meritorious performance of Spectral CUSUM through simulated and real-data examples of detecting changes in Yellowstone seismic sensors data.


The rest of the paper is organized as follows. Section \ref{sec:setup} provides a detailed formulation of the emerging communities problem as well as the switching membership problem. Section \ref{sec:CUSUM} introduces the exact-CUSUM and the proposed Spectral-CUSUM procedure. Section \ref{sec:Asymptotic} presents the asymptotic analysis of the proposed detection scheme together with parameter optimization and proof of the first-order asymptotic optimality. Section \ref{sec:subspace-track} presents an efficient gradient-based algorithm to keep track of the underlying community structure. Section \ref{sec:exp} gives simulation and real data examples to verify the theoretical findings and show the good performance of the proposed method. We delegate all proofs to the appendix.

\subsection{Literature}

{\it Community detection} in the {\it offline} setting (i.e., when the samples are collected beforehand and inference is made in one-shot) is a well-studied problem (see \cite{abbe2017community} for a survey). For example, spectral methods based on eigenvectors of the graph Laplacian are used in \cite{fiedler1973algebraic,pothen1990partitioning}. Besides, there are many practical algorithms for community detection (see, e.g., \cite{fortunato2010community}), including the so-called Kernighan–Lin algorithm \cite{kernighan1970efficient}, which uses a greedy algorithm to improve an initial division of the network and a genetic algorithm named Ga-net \cite{pizzuti2008ga} to detect the community structure by calculating the community score. However, offline community detection cannot be used to detect community changes if the graph is dynamic.

{\it Online community detection} has also been considered in the literature.
Peel and Clauset \cite{peel2015detecting} first formalized the change-point detection problem as identifying the times at which the large-scale patterns of interaction change fundamentally. They choose from a parametric family of probability distribution to describe the data and then use the Bayesian method to detect the change. A Markov-process-based approach is presented by \cite{wang2017fast}, which is based on MCMC. Each graph snapshot depends on the current generative model and the previously observed snapshot. Moreover, \cite{eswaran2018spotlight} proposes a method called Spotlight to detect anomalies in streaming graphs by composing a K-dimensional sketch containing K subgraphs to detect changes inside the dynamic graph. Recently, a Laplacian anomaly detection method for dynamic graphs is presented in \cite{huang2020laplacian}, which uses the spectrum of the Laplacian matrix of the graph structure at each snapshot to obtain low dimensional embeddings. However, none of these works gives a statistical perspective and asymptotic analysis for the quickest change detection in communities. While estimating the community structure through a dynamic network, our paper presents novel methods of handling estimation and change detection, which are supported by theoretical asymptotic optimality.


{\it Spectral graph change detection} is also related to our work. The spectral property of the graph is one of the essential theories that can capture the community structure of the graph. The spectral method proposed by  \cite{hewapathirana2020change} is used to detect changes in Noisy Dynamic Networks efficiently by transforming a graph to a lower-dimensional latent Euclidean space. In our work, we also use the spectral method to get a reduced dimensional representation for each node, and our procedure can detect different types of changes and has theoretically proved optimally.


{\it Change Detection of Gaussian Graphical Model} is similar to our work as well. We consider the Gaussian graphical model to capture the network structure through the inverse covariance matrix of the nodal features. The Gaussian model is useful for modeling the correlation between observations at different nodes. There are many previous works that are focused on the estimation of the Gaussian graphical model. For instance, penalized likelihood method is proposed by \cite{yuan2007model} for estimating the concentration matrix in the
Gaussian graphical model while \cite{liu2017tiger} presents an alternative tuning-insensitive approach to efficiently choose the tuning parameter in finite sample settings. But none of these estimation methods of the Gaussian graphical model can be directly used to detect the structure change. Recently, a piecewise stationary graphical model has been presented in \cite{keshavarz2020sequential}, and it is able to detect the change of the graph by monitoring the conditional log-likelihood of all nodes in the network. However, this method does not consider the community structure of the network and thus is unable to distinguish different types of changes inside the network.

{\it Manifold change detection} is another type of related work. In the work of \cite{xie2012change}, multi-scale online manifold learning is used to extract change-point detection test statistics from high-dimensional data. But they do not consider the network property for the high-dimensional data and thus cannot be applied directly to graph scenarios.

\section{Problem Setup}
\label{sec:setup}

Consider a dynamic network with $n$ nodes. Assume at each time $t$, we observe a feature value $v_{ti}$ for each node $i$, $i = 1, \ldots, n$. We collect all features at each time into one vector and denote as $v_t \in \mathbb R^n, t=1, 2, \ldots$. We will focus on the online setting where we observe feature values sequentially. Such a setting is widely applicable in real datasets, see below for some examples. 
\begin{itemize}
	\item In sensor networks, the feature $\{v_{ti},t=1,2,\ldots\}$ represent a sequence of signals recorded by the $i$-th sensor. The sensors may form communities, and thus the features can be correlated. It can be the observed seismic/solar activity measurement from each sensor at each time in a seismic/solar system. 
	
	\item In social networks, a vector of features represents user activities at each time. For instance, it can be social activities such as twittering at time $t$ from each user in the Twitter network.
\end{itemize}


Assume the features for nodes within the same community have a higher correlation than those that are not in the same community. The correlation can be estimated using observed nodal features \cite{yang2013community}. The underlying community structure may change at some time, which leads to a change in the correlations between the affected nodes and, thus the distribution of feature vectors. We aim to detect such a change as quickly as possible from sequential observations.

In the following, we first give the community definition in Section \ref{sec:community-def}  and statistical modeling for the community structure within graphs, then discuss two kinds of changes respectively in Section \ref{sec:community-change-def}.

\subsection{Adjacency matrix for community}
\label{sec:community-def}

Suppose there are $m$ communities within a network with $n$ nodes. Denote these communities as $m$ sets $\{C_1, \cdots, C_m\}$, where the $k$-th community is represented using the index set $C_k$ of the nodes belonging to this community. Assume the communities do not overlap with each other, i.e., the index sets are mutually disjoint. For each node $i$, we introduce an indicator vector $a_i \in\{0,1\}^{m}$ representing its true membership: the $k$-th entry equal to 1 and all other entries equal to 0 if node $i$ belongs to the $k$-th community $C_k$, i.e.,
$$
a_i = [0\ \cdots \! \underbrace{1}_{\text{$k$-th entry}}\! \cdots \ 0]^\top,\, \forall i\in C_k.
$$
Define the global indicator matrix as
\begin{equation}\label{eq:indicator}
A = \begin{bmatrix}a_1, \cdots ,a_n\end{bmatrix}^\top \in \{0,1\}^{n\times m}.
\end{equation}
Notice that we have
$$a_i^\top a_j = 
\begin{cases}
	1& \quad \text{if $\exists k$ s.t.  nodes}\  i,j \in C_k, \\
	0& \quad \text{otherwise.} \\
\end{cases}
$$
Therefore, the matrix $AA^\top\in\{0,1\}^{n\times n}$ defines an {\it adjacency matrix}, whose $(i,j)$-th entry equals to 1 if and only if nodes $i$ and $j$ belong to the same community. 
This setup for $A$ can also be generalized beyond 0-1 matrices. For example, we may let $a_i\in \mathbb{R}^{m}$ represent a feature embedding vector for node $i$, and each entry of $a_i$ represents the weight/probability for node $i$ belonging to the corresponding community.

\subsection{Community change-point detection}
\label{sec:community-change-def}

\begin{figure}[h!]
	\begin{center}
		\subfigure[Emergence problem illustration]{
			\includegraphics[width = 0.45\textwidth]{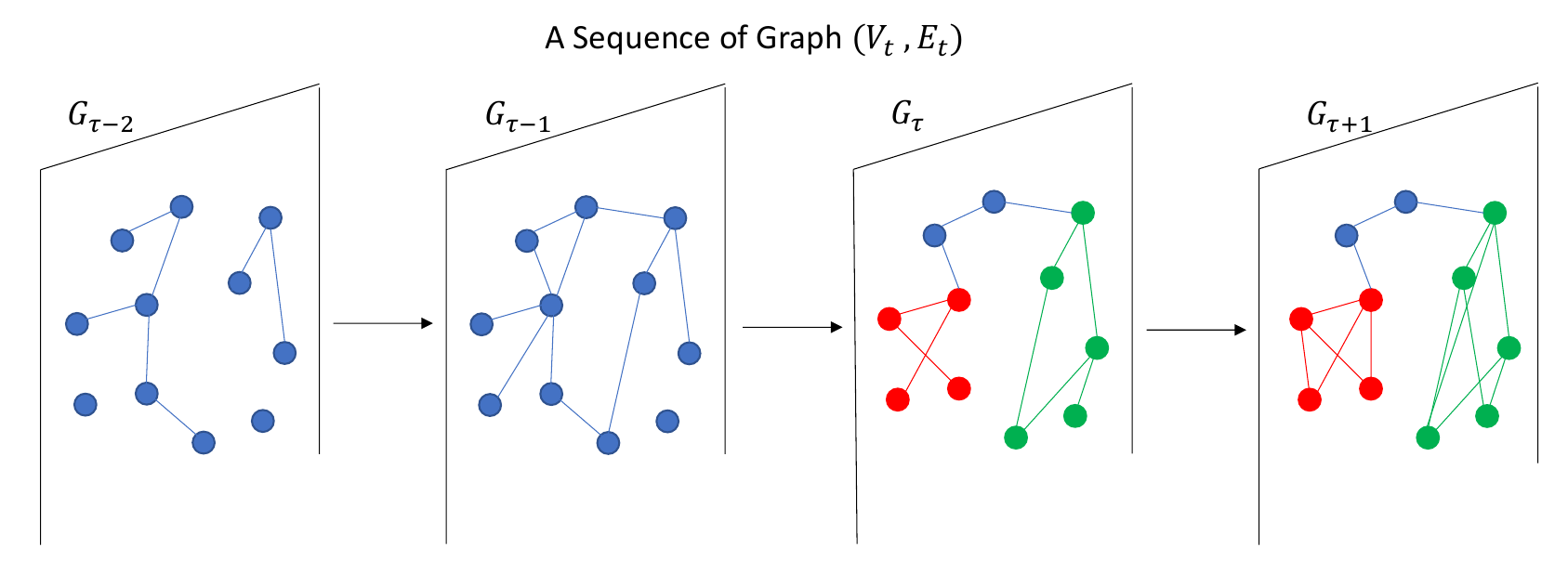}}
		\subfigure[Switching membership problem illustration]{
			\includegraphics[width = 0.45\textwidth]{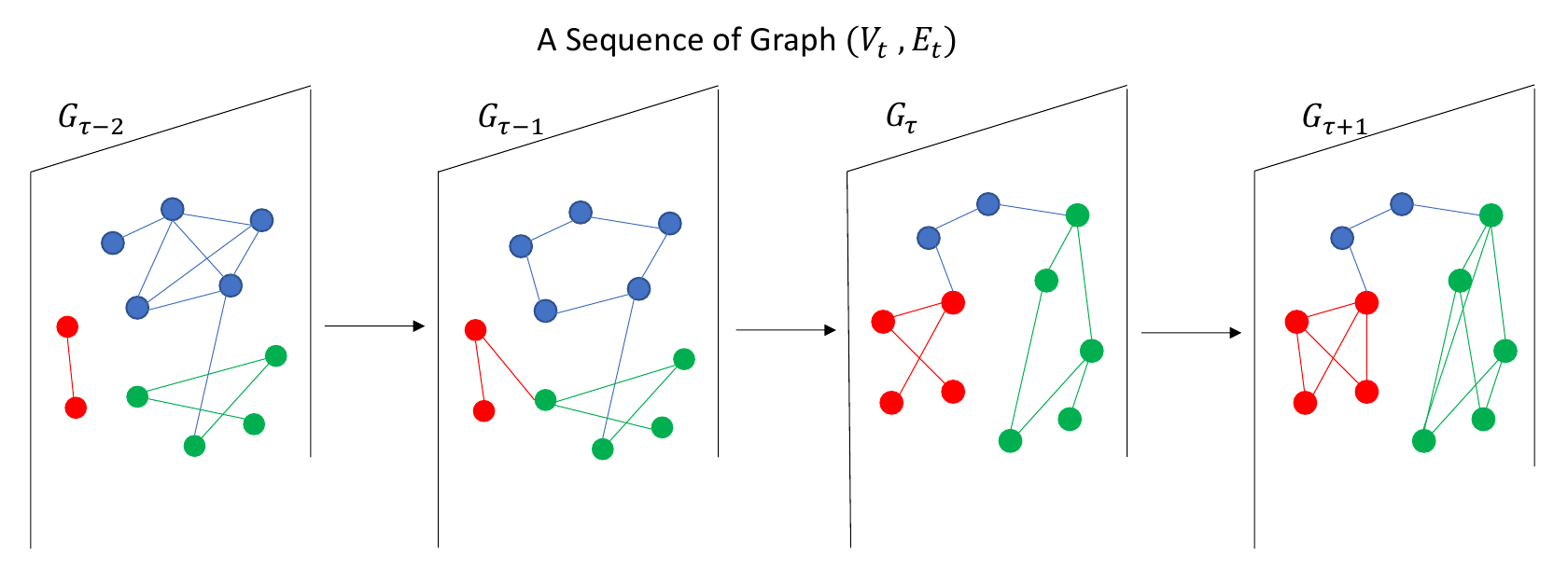}}
	\end{center}
	\caption{Case of the emergence of communities after change time $\tau$ as shown in (a) and switching membership case shown in (b). Three communities are marked as red, blue, and green dots. For the emergence problem, it can be seen that there is no community or one single community before changes happen, then after that, three separate communities emerge. For the switching membership case, it can be seen that before the change happens, the blue community is the biggest one, while after the change, the membership of it begins switching to the other two communities.}
	\label{fig:emerge}
\end{figure}

We aim to monitor two types of structural changes in dynamic networks:
\begin{enumerate}[label=\roman*)]
	\item The emergence of new communities: Before the change happens, there is no clear community formed in the graph, and the community structure emerges after the change, as indicated in Figure \ref{fig:emerge}(a).
	\item Switching membership: Some community members are switched after the change, such as the increase or decrease of a single community or the membership flow from one community to another, as indicated in Figure \ref{fig:emerge}(b).
\end{enumerate} 

We now formulate such change detection problems based on the adjacency matrix representation in Section \ref{sec:community-def}. 

\subsubsection{Emerging community}

In many applications, the change can be modeled as the emergence of several disjoint communities; the nodes inside the same community are more correlated with each other. Thus, we start by considering the emerging community detection problem, which assumes that the network has no community structure at the beginning but forms $m$ communities ($C_1, \cdots, C_m$) after the change where $C_k$ is the node index set for the $k$-th community. 

With the indicator matrix $A$ defined in \eqref{eq:indicator}, we assume the feature observations $v_t$ are multivariate Gaussian and its covariance matrix is modeled based on the adjacency matrix $AA^\top$. More specifically, we cast the emerging community problem as follows.
\begin{equation}
	\label{eq:model}
	\begin{array}{ll}
		H_0:  &v_{t} \distas{\text{i.i.d.}} \mathcal{N}(0, \frac{1}{\sigma^2}{I}), ~~~t=1,2,\dots\\
		H_1:&  v_{t} \distas{\text{i.i.d.}} \mathcal{N}(0,\frac{1}{\sigma^2} {I}), ~~~t=1,2,\dots,\tau,\\
		&v_{t} \distas{\text{i.i.d.}} \mathcal{N}(0,(AA^\top + \sigma^2 I)^{-1}),~~~ t = \tau+1,\tau+2,\dots\\
	\end{array}
\end{equation}
Here we introduce $\sigma^2 I$ as a noise term since the network is usually not perfectly separable in practice. Note that in the formulation here and below for the switching subspace, we consider the structure of the inverse covariance matrix, which is a common approach for Gaussian graphical models \cite{loh2013structure}. Following this, we can regard $AA^\top$ as the underlying community structure of the graph and $\sigma^2$ as the noise level. 

The above representation can be related to the Gaussian graphical model, which is a common approach to exploring the relationships between nodes in an undirected graph through the inverse covariance matrix. Given a Gaussian graphical model with covariance matrix $\Sigma$, there is an edge between node $i$ and node $j$, if and only if $\Sigma^{-1}_{ij}\neq 0$. Note that this can be related to our model: the zero-valued entry in the inverse covariance matrix means that the corresponding edge does not exist in the graphical model. As the inverse covariance matrix changes from $\sigma^2 I$ to $AA^\top + \sigma^2 I$ in a model \eqref{eq:model}, the corresponding off-diagonal entries change from zero to non-zero after the community emerges. 

\subsubsection{Switching membership} 

Another type of community change is called the switching membership problem. As shown in Figure \ref{fig:emerge}(b), some of the nodes belong to one community at first and switch to a different one after the change happens at time $\tau$. 
%
%
%
Similar to \eqref{eq:model}, the switching membership problem can be formulated as follows:
\begin{equation}
	\label{eq:model2}
	\begin{array}{ll}
		H_0:  &v_t \distas{\text{i.i.d.}} \mathcal{N}(0,(A_1A_1^\top + \sigma^2 I)^{-1}), ~~~t=1,2,\dots\\
		H_1:&  v_t \distas{\text{i.i.d.}} \mathcal{N}(0,(A_1A_1^\top + \sigma^2 I)^{-1}), ~~~t=1,2,\dots,\tau,\\
		&v_t \distas{\text{i.i.d.}} \mathcal{N}(0,(A_2A_2^\top + \sigma^2 I)^{-1}),~~~ t = \tau+1,\tau+2,\dots\\
	\end{array}
\end{equation}
Here $A_1$ represents the pre-change community structure while $A_2$ represents the post-change community structure. This general model can denote cases when the sizes of certain communities change. It can also model the case when the total number of communities increases (one community splits into smaller ones) or decreases (several small communities merge into a bigger one). As a result, the emergence problem can be seen as a special case of the switching membership problem, which is capable of detecting various types of graph structure changes. We will show later the detecting procedure of switching membership can also be equivalently treated as a community emergence problem.

%

\section{Detection Procedures} \label{sec:CUSUM}

In this section, we first review the well-known cumulative sum (CUSUM) detection rule and then propose the Spectral-CUSUM procedure under both emergence and switching membership scenarios. 

\subsection{Exact-CUSUM procedure}

Let $f_{\infty}(\cdot)$ and $f_0(\cdot)$ denote the pre- and post-change probability density function (pdf) of the observations, and $\mathbb{E}_\infty$ and $\mathbb{E}_0$ denotes the expectation under $f_\infty$ and $f_0$, respectively. The CUSUM statistic \cite{page1954continuous} is defined by maximizing the log-likelihood ratio statistic over all possible change-point locations:
\begin{equation*}
	S_t = \max_{1\leq k\leq t}\sum_{i=k}^t\log \frac{f_0(v_i)}{f_{\infty}(v_i)}.
\end{equation*}
$S_t$ has a recursive formulation with $S_0 = 0$ as follows:
\begin{equation*}
	\label{eq5}
	S_t = (S_{t-1})^+ +\log \frac{f_0(v_t)}{f_{\infty}(v_t)}, \ t\geq 1,
\end{equation*}
where $(x)^+:=\max\{x,0\}$. The corresponding CUSUM stopping time $T_C$ is defined as:
\begin{equation}\label{eq:cusumstop}
	T_C = \inf\{t>0:S_t\geq b\},
\end{equation}
where $b$ is a pre-set constant threshold.
Under the model (\ref{eq:model}), we have that
\begin{equation*}
	\label{eq6}
	\begin{aligned}
		\log \frac{f_0(v_t)}{f_{\infty}(v_t)} &=-\frac{1}{2}v_t^\top AA^\top v_t+\frac{1}{2}\log\frac{\mbox{det}(AA^\top+\sigma^2 I)}{\sigma^{2n}}.\\
	\end{aligned}
\end{equation*}

Since the multiplicative factor $1/2$ is positive, we can omit it from the log-likelihood ratio when forming the CUSUM statistic, thus yielding an equivalent formulation:
\begin{equation}
	\label{eq:CUSUM}
	S_t = (S_{t-1})^+ -v_t^\top AA^\top v_t + \underbrace{\log \frac{\mbox{det}(AA^\top+\sigma^2 I)}{\sigma^{2n}}}_{d},
\end{equation}
where $d$ is a drift parameter that is fixed in the CUSUM procedure. 

Similarly, for the switching membership problem under the model (\ref{eq:model2}), the log-likelihood ratio is:
\begin{equation*}
	\label{eq6-2}
	\begin{aligned}
		\log \frac{f_0(v_t)}{f_{\infty}(v_t)} 
		=&-\frac{1}{2}v_t^\top(A_2A_2^\top-A_1A_1^\top)v_t\\
		&+\frac{1}{2}\log\frac{\mbox{det}(A_2A_2^\top+\sigma^2 I)}{\mbox{det}(A_1A_1^\top+\sigma^2 I)}.\\
	\end{aligned}
\end{equation*}
Therefore, the CUSUM statistic in this case becomes:
\begin{equation}
	\label{eq:CUSUM-2}
	S_t = (S_{t-1})^+ -v_t^\top(A_2A_2^\top-A_1A_1^\top)v_t +\underbrace{\log \frac{\mbox{det}(A_2A_2^\top+\sigma^2 I)}{\mbox{det}(A_1A_1^\top+\sigma^2 I)}}_{d^\prime}.
\end{equation}
We note that in the exact CUSUM procedure, all the parameters are assumed known so that the drift terms $d$ and $d'$ can be computed explicitly beforehand.  

\subsection{Spectral-CUSUM procedure}

The implementation of the exact CUSUM procedure requires that all parameters are known, and it has been proved to be optimum \cite{lorden1971procedures,moustakides1986optimal}. However, if the post-change distribution is unknown, the exact-CUSUM is not applicable. Usually, we can estimate the pre-change distribution using historical data (training data), but the post-change community structures are unknown since it represents anomaly information and cannot be predicted. Therefore, the post-change distribution, i.e., the post-change structure $A$, has to be estimated sequentially from data. A natural estimate for the post-change covariance matrix $(AA^\top + \sigma^2 I)^{-1}$ is the sample covariance matrix. To eliminate the matrix inversion computation in estimating $A$, we use eigendecomposition on the sample covariance matrix to estimate $A$ directly. 
Therefore, we propose the Spectral-CUSUM procedure below, where we approximate the covariance matrix by rank-$m$ eigendecomposition.


Define $U\in\mathbb{R}^{n\times m}$ and $\Lambda\in\mathbb{R}^{m\times m}$ as the eigenvectors and eigenvalues of the term $A A^\top$ within the post-change covariance matrix, i.e., $A A^\top=U\Lambda U^\top$. Here $\Lambda$ is a diagonal matrix corresponding to the eigenvalues of the matrix $A A^\top$ defined using \eqref{eq:indicator}:
\[
\Lambda =
\begin{bmatrix}
	|C_{1}| & & \\
	& \ddots & \\
	& & |C_{m}|
\end{bmatrix},
\]
where $|C_k|$ denotes the number of nodes inside the $k$-th community.
This shows that the eigenvalue decomposition of the adjacency matrix would reflect the community structure of a graph, which can also provide empirical guidance on how to determine the potential number $m$ of communities.


We can then jointly estimate matrices $U$ and $\Lambda$. Using observations $\{v_{t+1},\dots,v_{t+w}\}$ in a {\it future} sliding window with length $w$, we define
\begin{equation}
	\label{eq:emp_est}
	\hatG_{t}=(v_{t+1}v^\top_{t+1}+\dots+v_{t+w} v_{t+w}^\top)/w.
\end{equation}
Note that $\hatG_t$ serves as an approximation for the covariance matrix $(AA^\top+\sigma^2I)^{-1}$ if samples $\{v_{t+1},\ldots,v_{t+w}\}$ are drawn from post-change distribution; for ease of presentation, we assume eigenvalues are distinct. Let $\{\hat{u}_{t1},\dots,\hat{u}_{tm}\}$ be the $m$ unit-norm eigenvectors corresponding to the $m$ {\it smallest} eigenvalues $\{\hat{\lambda}_{t1}<\dots<\hat{\lambda}_{tm}\}$ of $\hatG_t$. Let
$$\hatU_t = [\hat{u}_{t1}, \dots, \hat{u}_{tm}] \in \mathbb{R}^{n\times m},$$ and
$$\hatLambda_t = \text{diag}(\hat{\lambda}_{t1},\dots,\hat{\lambda}_{tm}) \in \mathbb{R}^{m\times m}.$$
Then $A$ can be approximately estimated using $\widehat{A}_t=\hatU_t\hatLambda_t^{-1/2}$, as the noise term $\sigma^2$ is usually relatively small. Then we can substitute the estimate $\widehat{A}_t$ into \eqref{eq:CUSUM} to obtain alternative detection statistics, which we call Spectral-CUSUM, when the post-change distribution is unknown:
\begin{equation}
	\label{eq:subspaceCUSUM}
	\mathcal{S}_t = (\mathcal{S}_{t-1})^+ - v_t^\top \widehat A_t\widehat A_t^\top v_t+ d.
\end{equation}
Here $d$ is a tunable drift parameter that plays a similar role as the last term in \eqref{eq:CUSUM}. And the stopping time is defined as follows 
\begin{equation}
	\label{eq:CUSUMstop}
	\mathcal{T}_C = \inf\{t>0:\mathcal{S}_t\geq b\}.
\end{equation}
Similarly, for the switching membership problem, the alternative to \eqref{eq:CUSUM-2} is: 
\begin{equation}
	\label{eq:switchingCUSUM}
	\mathcal{S}_t = (\mathcal{S}_{t-1})^+ -v_t^\top(\widehat A_t \widehat A_t^\top-A_1A_1^\top)v_t+d^\prime.
\end{equation}




In the proposed Spectral-CUSUM procedure, the drift parameter $d$ should be chosen to ensure the detection statistics are capable of detecting the change. More specifically, the CUSUM type procedure requires the increment term in $\mathcal{S}_t$ to have a {\it negative} mean under the pre-change distribution, and a {\it positive} mean under the post-change distribution. Therefore, for the emerging community problem, we need:
\begin{equation}
	\label{eq:CUSUMrequire}
	\mathbb{E}_0[v_t^\top \widehat A_t \widehat A_t^\top v_t]<d<\mathbb{E}_\infty[v_t^\top \widehat A_t \widehat A_t^\top v_t].
\end{equation}
For the switching membership problem, the same rule means that we need:
\begin{equation}
	\label{eq:CUSUMrequire2}
	\mathbb{E}_0[v_t^\top (\widehat A_t \widehat A_t^\top-A_1A_1^\top) v_t]<d'<\mathbb{E}_\infty[v_t^\top(\widehat A_t \widehat A_t^\top-A_1A_1^\top)v_t].
\end{equation}
Since $A_1$ is a constant matrix known beforehand from historical data, we can calculate the expectations  $\mathbb{E}[v_t^\top A_1A_1^\top v_t]$ explicitly. Thus, the switching membership problem can be treated as an emerging community problem. Therefore, the theoretical analysis for the switching membership problem is similar to the emerging community problem, and we will only discuss the emergence problem in Section \ref{sec:Asymptotic}. 

Due to the above property \eqref{eq:CUSUMrequire}, which is also mentioned in \cite{oskiper2005quickest,egea2017comprehensive}, the detection statistic in \eqref{eq:subspaceCUSUM} will deviate from 0 and increase gradually after the change happens.

\begin{algorithm}[H]
	\caption{Spectral-CUSUM procedure}\label{Spectral-CUSUM-alg}
	\begin{algorithmic}[1]
		\Require{ Sequence of observations $\{v_t, t=1,2,\dots\}$, number of communities $m$, sliding window size $w$, carefully selected drift parameter $d$ (or $d'$), detection threshold $b$.}
		\Ensure{Stopping time $\mathcal{T}_C$.}
		\State Initialize $\mathcal{T}_C=\infty$, $t=0$, $S_0=0$.
		\While{$S_t < b $}
		\State Calculate sample covariance matrix $\hatG_t$ from future observations using (\ref{eq:emp_est});
		\State Compute $\hatU_t$ and $\hatLambda_t$ through the eigenvalue decomposition of $\hatG_t$;
		\State Estimate community structure $\widehat{A}_t=\hatU_t\hatLambda_t^{-1/2}$;
		\State Let $t=t+1$ and update Spectral-CUSUM statistics $\mathcal{S}_t$ following (\ref{eq:subspaceCUSUM}) or (\ref{eq:switchingCUSUM}).
		
		\EndWhile
		\State Set $\mathcal{T}_C=t$.
		\State \Return {$\mathcal{T}_C$}
		
	\end{algorithmic}
\end{algorithm}

\section{Theoretical Analysis}
\label{sec:Asymptotic}

This section provides a theoretical analysis of the proposed Spectral-CUSUM procedure under the emerging community setting. The main result is presented in Theorem \ref{thm:main}, which shows the asymptotic optimality of Spectral-CUSUM under the optimal choices of parameters. We also derive the form of the optimal parameters. The analysis techniques are related and extended to those used in prior works in \cite{XieMoustakidesXie2022} and in \cite{xie2018first,xie2020sequential} for subspace change detection.

\subsection{Preliminary}

We first state the following assumptions that are being made in order to derive the main results. 
\begin{assumption}
	\label{assumption}
	We make the following assumptions for the true community structure.
	\begin{enumerate}
		\item[(1)] The total number of post-change communities remains a constant $m$. 
		\item[(2)] 
		\textcolor{black}{(Community sizes are comparable but not
		not identical) The size of different communities are not identical and is ordered as $|C_1|>|C_2|>\cdots>|C_m|$, and $|C_1|/|C_m| \leq (1+\theta)^2$ for some $\theta > 0$.}
		\item[(3)] \textcolor{black}{(Noise variance bounded) 
$$\frac{\sigma^2}{|C_i|} <  
\left(\frac 1 4 +
\eta \cdot \sqrt{\frac{{m-1}}{{n-m}}}\cdot \frac{1}{\left(
    1+\theta - \frac{1}{1+\theta}
    \right)^2}
\right)^{1/2} - \frac{1}{2}$$, $\forall i$, for a small constant $0<\eta < 1$}.
	\end{enumerate}
\end{assumption}
\textcolor{black}{We remark that Assumption \ref{assumption}(1) can be treated as a pre-defined number of communities to detect. Assumption \ref{assumption}(2) assumes that the community sizes are comparable. Assumption \ref{assumption}(3) requires that the noise variance relative to the community sizes is not too large, but we do not need the noise variance to be 0. In fact, when $\theta \rightarrow 0$, i.e., when the community size is very similar to each other, the ratio $\frac{1}{(1+\theta - \frac{1}{1+\theta})^2} \rightarrow \infty$ and thus the condition can be easily satisfied; we can even allow high noise $\sigma^2$ when the community sizes are comparable.}


Define the following quantities to simplify the presentation of the main results:
\begin{equation}
\begin{aligned}
	B_i&=\sum_{k=1,k\neq i}^{m}\frac{\lambda_i \lambda_k}{(\lambda_i-\lambda_k)^2} >0,\\
		D &= \sum_{i=1}^m \frac{\lambda_i}{\lambda_i+\sigma^2} \bigg(1-\frac{B_i}{w}+\frac{3B_i^2}{w^2}\bigg),\\
	\widetilde D &= \sum_{i=1}^m \frac{\lambda_i}{\lambda_i+\sigma^2} \bigg(1-\frac{B_i}{w}\bigg);\\
	\label{BD_def}
\end{aligned}
\end{equation}
recall that $\lambda_i=|C_i|$ is the $i$-th largest eigenvalue of the matrix $AA^\top$ (thus we always have $\lambda_i>0$), 
$w$ is the sliding window size, and $\sigma^2$ is the noise level. For theoretical analysis purposes, we assume the eigenvalues are distinct, i.e., no community sizes are exactly the same (see Assumption \ref{assumption}(2)); not that when this does not hold, we can still apply the algorithm although the theory needs to be extended. \textcolor{black}{Because of Assumption \ref{assumption}(2), we have 
\begin{equation}
\textcolor{black}{
\frac{m-1}{ \Big(1+\theta - \frac1{1+\theta} \Big)^2 } <B_i<\frac{m-1}{ \Big(\sqrt{\frac{|C_1|+1}{|C_1|}} - \sqrt{\frac{|C_1|}{|C_1|+1}} \Big)^2 }, \ \forall i=1,\ldots,m.}
\label{B_bound}
\end{equation}
Moreover, a basic bound can be derived: $\textcolor{black}{D \geq \frac{11}{12} \sum_{i=1}^m \frac{\lambda_i}{\lambda_i+\sigma^2} > 0}$ based the property of the quadratic equation, \textcolor{black}{and $\frac{11}{12} \sum_{i=1}^m \frac{\lambda_i}{\lambda_i+\sigma^2}\leq \textcolor{black}{\widetilde D} \leq D\leq \sum_{i=1}^m \frac{\lambda_i}{\lambda_i+\sigma^2}<m$ when $w$ is sufficiently large}.
}

We start with a useful result from \cite{anderson1963asymptotic}.
\begin{theorem}[Asymptotic Property of Sample Covariance, \cite{anderson1963asymptotic}]
	\label{thm1}
	The asymptotic distribution of eigenvectors and eigenvalues of sample covariance matrix \eqref{eq:emp_est} from $w$ samples under post-change distribution \eqref{eq:model} are: (i) independent and (ii) has the following distribution under the assumption that the eigenvalues are distinct: 
	\begin{equation*}
		\label{eq:asy_eigenvector}
		\begin{aligned}
			\sqrt{w}(\hat{u}_{i}-u_i)
			&\stackrel{d}{\longrightarrow}\mathcal{N}\bigg(0, \sum_{k=1,k\neq i}^{m}\frac{\lambda_i \lambda_k}{(\lambda_i-\lambda_k)^2}u_ku_k^\top\bigg),
		\end{aligned}
	\end{equation*}
	\begin{equation*}
		\label{eq:asy_eigenvalue}
		\begin{aligned}
			\sqrt{w}(\hat{\lambda}_{i}-\lambda_i)\stackrel{d}{\longrightarrow}\mathcal{N}(0, 2\lambda_i^2),
		\end{aligned}
	\end{equation*}
	where $\stackrel{d}{\longrightarrow}$ denotes convergence in distribution \textcolor{black}{when $w\rightarrow \infty$}. Here $\hat{u}_{i}$ and $\hat{\lambda}_{i}$ denotes the sample eigenvectors and eigenvalues, while $u_{i}$ and $\lambda_{i}$ denotes the eigenvectors and eigenvalues of the true covariance matrix, and $n$ is the data dimension.
\end{theorem}

For the emerging community problem in \eqref{eq:model}, the post-change covariance matrix is $(AA^\top + \sigma^2 I)^{-1}$. Recall that we use the window size $w$ to construct the post-change sample covariance matrix in \eqref{eq:emp_est} and then estimate the eigenvalues and eigenvectors. Denote $\rho_1<\cdots<\rho_m$ as the smallest $m$ eigenvalues of $(AA^\top + \sigma^2 I)^{-1}$; based on the eigendecomposition $AA^\top = U\Lambda U^\top$, we can show $\rho_i=1/(\sigma^2 + |C_i|)$. Let $u_i$, $i=1,\ldots,m$ denote the corresponding eigenvectors. Let $\hat\rho_i$ and $\hat u_i$, $i=1,\ldots,m$, be the smallest $m$ eigenvalues and corresponding eigenvectors obtained from the sample covariance matrix $\hatG_t$. By Theorem \ref{thm1}, we have
\[
\sqrt{w}(\hat{\rho}_{i}-\rho_i)\stackrel{d}{\longrightarrow}\mathcal{N}(0, 2\rho_i^2),
\]
and
\begin{equation*}
\begin{aligned}
	\sqrt{w}(\hat{u}_{i}-u_i)
	\stackrel{d}{\longrightarrow}\mathcal{N}\bigg(0, &\sum_{k=1,k\neq i}^{m}\frac{\rho_i \rho_k}{(\rho_i-\rho_k)^2}u_ku_k^\top \\
	&+ \frac{\rho_i/\sigma^2}{(\rho_i-1/\sigma^2)^2} (I- UU^\top)\bigg).
	\label{u_est_error}
\end{aligned}
\end{equation*}
\textcolor{black}{We can show in the above estimation error covariance matrix \eqref{u_est_error}, the second term in the summation is negligible compared with the first term under Assumption\ref{assumption}, when $\eta$ is sufficiently small (e.g., 0.1, 0.2) --- $\eta$ is desired upper bound for the ratio of the second term relative to the first term. The derivation can be found in the appendix.
} 

%

With the help of Theorem \ref{thm1} we can show the following properties for Spectral-CUSUM.
\begin{lemma}[Properties of detection statistic for Spectral-CUSUM]
	\label{lemma:1}
	\textcolor{black}{When $w\rightarrow \infty$,} the expected drifts under the pre- and post-change distributions for Spectral-CUSUM are given by
	\begin{equation*}
\mathbb{E}_{\infty}[v_t^\top \widehat{A}_t\widehat{A}_t^\top v_t] 
			=m, \textcolor{black}{\quad \widetilde D\leq\mathbb{E}_{0}[v_t^\top \widehat{A}_t\widehat{A}_t^\top v_t]\leq D}.
	\end{equation*}
\end{lemma}
The proof of this Lemma is shown in the Appendix. 
\textcolor{black}{By Lemma\,\ref{lemma:1}, the necessary condition for the drift parameter $d$ in \eqref{eq:CUSUMrequire} translates into $D<d<m$. Later in Section \ref{sec:opt}, we can show that by construction $\mathbb{E}_{\infty}[v_t^\top \widehat{A}_t\widehat{A}_t^\top v_t] > \mathbb{E}_{0}[v_t^\top \widehat{A}_t\widehat{A}_t^\top v_t]$ (i.e., $m>D$) can be satisfied, and hence we can choose a suitable drift parameter $d$ in between such that the Spectral-CUSUM procedure work. \textcolor{black}{Indeed, we emphasize that we focus on the regime that $w$ is sufficiently large such that $D<m$ is guaranteed.}}

%
%


\subsection{ARL/EDD performance analysis}

The standard performance of change detection procedures is measured by average run length (ARL) and expected detection delay (EDD). ARL represents the average time interval between two consecutive false alarms, while EDD measures the worst-case detection delay. When ARL is fixed, it is known that the exact CUSUM procedure minimizes EDD, which can be calculated directly. \textcolor{black}{In the following, we will analyze ARL and EDD for the proposed Spectral-CUSUM procedure given in \eqref{eq:subspaceCUSUM} for the emerging community problem, under the assumption that the window size $w^*$ scales in order $\sqrt{\log\gamma}$ in the asymptotic analysis as $\gamma\rightarrow \infty$.} \textcolor{black}{In this case, $D/\widetilde D$ approaches to 1, and then we can approximate $\mathbb{E}_{0}[v_t^\top \widehat{A}_t\widehat{A}_t^\top v_t]$ by $D$ in our derivation. }

Given a constant $\gamma > 1$ as the desired lower bound of ARL, we need to set the threshold $b$ in \eqref{eq:cusumstop} (and similarly in \eqref{eq:CUSUMstop} for Spectral-CUSUM) accordingly such that $\text{ARL} \geq \gamma$. 
Recall that $T_C$ denotes the stopping time for the exact CUSUM procedure. Thus $\mathbb{E}_0[T_C]$ and $\mathbb{E}_\infty[T_C]$ are EDD and ARL of CUSUM. Similarly, $\mathbb{E}_0[\mathcal{T}_C]$ and $\mathbb{E}_\infty[\mathcal{T}_C]$ are the EDD and ARL of Spectral-CUSUM. According to classic results from \cite{siegmund2013sequential}, we have the following:
\begin{equation}\label{eq:thmcusum}
	\mathbb{E}_\infty[T_C]=\frac{e^b}{\mathcal{I}_\infty}(1+o(1)),\   \mathbb{E}_0[T_C]=\frac{b}{\mathcal{I}_0}(1+o(1)), 
\end{equation} 
where $\mathcal{I}_0$ and $\mathcal{I}_\infty$ are the Kullback-Leibler (KL) divergence: %
\[
\mathcal{I}_\infty= \mathbb{E}_\infty\{\log[f_\infty(x)/f_{0}(x)]\}, \quad 
\mathcal{I}_0 = \mathbb{E}_0\{\log[f_0(x)/f_{\infty}(x)]\},
\]
The constraint $\text{ARL}\geq\gamma$ will be satisfied with threshold $b = (\log \gamma)(1 + o(1))$ according to \eqref{eq:thmcusum}.
\textcolor{black}{
\begin{lemma}[KL divergence for emerging subspace case]\label{lemma:KL}
For emerging subspace problem \eqref{eq:model}, we have the Kullback–Leibler divergence (K-L) divergence
\[
    \mathcal{I}_0 = \mathbb{E}_0\bigg[\log \frac{f_0(v)}{f_{\infty}(v)}\bigg] = 
    -\frac{1}{2}\sum_{i=1}^m h\left(
    \frac{\lambda_i}{\sigma^2+\lambda_i} 
    \right),
\]
where $h(x) = x+\log(1-x)$.
\end{lemma}
Note that $h(x) < 0$, for $x\in (0, 1)$, so this also verifies that $\mathcal I_0 >0$.
}
Using the K-L divergence in Lemma \ref{lemma:KL}, we can obtain the EDD expression:
\begin{equation}
\label{exact_CUSUM}
\begin{split}
	\mathbb{E}_0[T_C]&=\frac{\log\gamma}{\mathcal{I}_0}(1+o(1)).
	\end{split}
\end{equation}

%
For performance analysis of Spectral-CUSUM, we follow a similar strategy as \cite{xie2018first}, where the analysis is done for a different subspace detection problem under the special rank-one case; here, to generalize the analysis, we extend it for rank more than one. Since the increment term $-v_t^\top\hatU_t\hatLambda^{-1}_t\hatU_t^\top v_t + d$ in (\ref{eq:subspaceCUSUM}) is not a log-likelihood ratio. Thus we cannot use the ARL and EDD expressions in 
\eqref{eq:thmcusum} directly, which are derived for log likelihood-ratio based CUSUM procedures. To compute the ARL and EDD of $\mathcal T_C$, we introduce an {\it equalizer} $\delta_\infty \textcolor{black}{\in \mathbb{R}}$ such that:
\begin{equation}
	\label{eq:equalizer}
	\mathbb{E}_\infty\left[\exp\{\delta_\infty(-v_t^\top \hatU_t\hatLambda^{-1}_t\hatU_t^\top v_t+ d)\}\right]=1.
\end{equation}
And then when \eqref{eq:equalizer} holds, $\delta_\infty[-v_t^\top \hatU_t\hatLambda^{-1}_t\hatU_t^\top v_t + d]$ is the log-likelihood ratio between $\tilde{f}_0$ and $f_\infty$ where 
\[\tilde{f}_0 = \exp\{\delta_\infty[-v_t^\top \hatU_t\hatLambda^{-1}_t\hatU_t^\top v_t+ d]\}f_\infty.\]
This allows us to compute the threshold $b$ asymptotically as $b = (\log \gamma)(1+o(1))/\delta_\infty$. 
Similarly, we can find a $\delta_0 > 0$ and define $\tilde f_\infty = \exp\{\delta_0[v_t^\top \hatU_t\hatLambda^{-1}_t\hatU_t^\top v_t- d]\}f_0$ so that $\delta_0[v_t^\top \hatU_t\hatLambda^{-1}_t\hatU_t^\top v_t- d]$ is the log-likelihood ratio between $f_0$ and $\tilde f_\infty$, leading to $\mathbb E_0[\mathcal T_C] = b(1+o(1))/(\mathbb{E}_0[-v_t^\top \hatU_t\hatLambda^{-1}_t\hatU_t^\top v_t]+d)$ where the dependence on $\delta_0$ being an $o(1)$ term. Now after substituting $b$, using the ARL and EDD expression for the likelihood-ratio based CUSUM in \eqref{eq:thmcusum}, we obtain:
\begin{equation}\label{eq:edd_1}
\begin{aligned}
\mathbb{E}_0(\mathcal{T}_C)&=\frac{\log\gamma\big(1+o(1)\big)}{\delta_\infty(\mathbb{E}_0[-v_t^\top \hatU_t\hatLambda^{-1}_t\hatU_t^\top v_t]+d)}+w\\ &= \frac{\log\gamma\big(1+o(1)\big)}{\delta_\infty(-D+d)}+w ,
\end{aligned}
\end{equation}
where the second equality is due to Lemma \ref{lemma:1}; and the window length $w$ is added for the reason that we are using additional data to perform estimation to detect the potential change at time $t$; thus, the actual detection time is $t+w$.

Now would like to find the equalizer $\delta_\infty$ that satisfies \eqref{eq:equalizer} and thus can lead to the desired likelihood ratio construction above. Using standard computations involving Moment Generating Function for Gaussian random variables, we can write:
\begin{equation}
	\label{eq:relat}
	\begin{split}
		&\mathbb{E}_\infty[e^{\delta_\infty[-v_t^\top \hat{A}_t\hat{A}_t^\top v_t + d]}]\\
		=&e^{\delta_\infty d}\mathbb{E}\bigg[\mathbb{E}_\infty[e^{-\delta_\infty[v_t^\top \hat{A}_t\hat {A}_t^\top v_t]}|\hat{A}_t]\bigg]\\
		=&e^{\delta_\infty d}\mathbb{E}\Bigg[\int e^{-\delta_{\infty}(v_t^\top \hat{A}_t\hat {A}_t^\top v_t)}\frac{e^{-v_t^\top v_t \sigma^2/2}}{\sqrt{(2\pi)^n(1/\sigma^2)^n }} dv_t\bigg]\\
		=&e^{\delta_\infty d}\mathbb{E}_\infty\left[ |I+2\sigma^{-2}\delta_\infty \hat{A}_t\hat {A}_t^\top|^{-1/2}\right]\\
		=&\frac{e^{\delta_\infty d}}{\sqrt{\prod_{i=1}^m(1 + 2\sigma^{-2}\delta_\infty/\rho_i)}}
		= 1.
	\end{split}
\end{equation}
In evaluating the integral above, we use the standard technique of
``completing the square" in the exponent for Gaussian distribution. With proper normalization, we generate an alternative Gaussian probability density function which integrates to 1. 
Solving (\ref{eq:relat}) we obtain that the drift term $d$ is related to the equalizer $\delta_\infty$:
$
d = (2\delta_{\infty})^{-1}\sum_{i=1}^m\log({1+2\sigma^{-2}\delta_\infty/\rho_i}).
$
We notice that under the pre-change measure, all eigenvalues $\rho_i$ of the covariance matrix are equal to $1/\sigma^2$, thus the drift parameter can be further written as:
\begin{equation}
	\label{eq:d_exp}
	d = \frac{m\log({1+2\delta_\infty})}{2\delta_{\infty}}.
\end{equation}
Using \eqref{eq:d_exp}, we can eliminate $d$ from the EDD expression and only leave $\delta_\infty$ dependence. \textcolor{black}{We will later validate that the drift $d$ in \eqref{eq:d_exp} under the optimal choice $\delta_\infty^*$ will indeed satisfy the  condition \eqref{eq:CUSUMrequire} to ensure a valid CUSUM procedure.} Combining Equation \eqref{eq:edd_1} and Equation \eqref{eq:d_exp} we have 
the expression for EDD is: 
\begin{equation}
	\label{eq:final_EDD}
	\begin{aligned}
		&\mathbb{E}_0(\mathcal{T}_C)=
		\frac{2\log\gamma\big(1+o(1)\big)}{-2\delta_\infty D + m\log ( 1+2\delta_\infty)}+w.
	\end{aligned}
\end{equation}
In the following, we will further derive the optimal value of $\delta_\infty$ (thus the optimal value of the drift parameter since they are equivalent through \eqref{eq:d_exp}) to minimize the EDD, as a function of $w$. 

\subsection{Optimal parameters for Spectral-CUSUM to minimize EDD}\label{sec:opt}

Note that the formulation \eqref{eq:final_EDD} contains two parameters: the drift $d$ (or the equalizer $\delta_\infty$) and the window size $w$. We will first optimize over these  parameters to minimize EDD, then show that the Spectral-CUSUM procedure based on optimal window size and draft is first-order asymptotically optimum.


%
We first find the optimal value of $\delta_\infty$ and the corresponding drift parameter $d$ (according to \eqref{eq:d_exp}). We observe that the denominator in (\ref{eq:final_EDD}) is a concave function of $\delta_\infty$ therefore, it exhibits a single maximum. Setting the derivative of the denominator as a function of $\delta_\infty$ to be 0, we obtain the optimum value of $\delta_\infty$ (we omit the high order terms of eigenvalue in the product):
\begin{equation}
	\label{optimal_delta}
	\delta_\infty^* 
	= \frac{m}{2D}-\frac{1}{2}.
\end{equation}
\textcolor{black}{Note that the optimal drift $d^* = \frac{m\log({1+2\delta_\infty}^*)}{2\delta_{\infty}^*}$ corresponding to the optimal $\delta_\infty^*$ satisfies the condition \eqref{eq:CUSUMrequire}. 
Substitute the optimal $\delta_\infty^*$, we can show that 
\begin{equation}
D< \frac{m\log(1+2\delta_\infty^*)}{2\delta_\infty^*} = \frac{m\log(\frac{m}{D})}{\frac{m}{D}-1} < m,\label{D_cond}
\end{equation}
where we have used that $x\log x \geq x-1$ for $x>1$. By letting $x=m/D$, we have $x>1$ and prove the left hand side of the inequality; similarly by $\log(x)<x-1$ for $x>1$, we have the right hand side of the inequality.}

Substituting $\delta_\infty^*$ to \eqref{eq:final_EDD}, we obtain the optimized EDD (with respect to $\delta_\infty$) as a function of $w$:
\begin{equation}
	\label{eq:EDD_min}
	\mathbb{E}_0(\mathcal{T}_C)=\frac{2\log\gamma\big(1+o(1)\big)}{m g(D/m)}+w.
\end{equation}
 \textcolor{black}{where $g(x) =x-1-\log(x)$. From \eqref{D_cond} we know $D/m <1$. 
It should be noted that in the expression above, the dependence on $w$ also comes from $D$ defined in \eqref{BD_def}. }


\begin{figure}
    \begin{center}
      \includegraphics[width=0.5\textwidth]{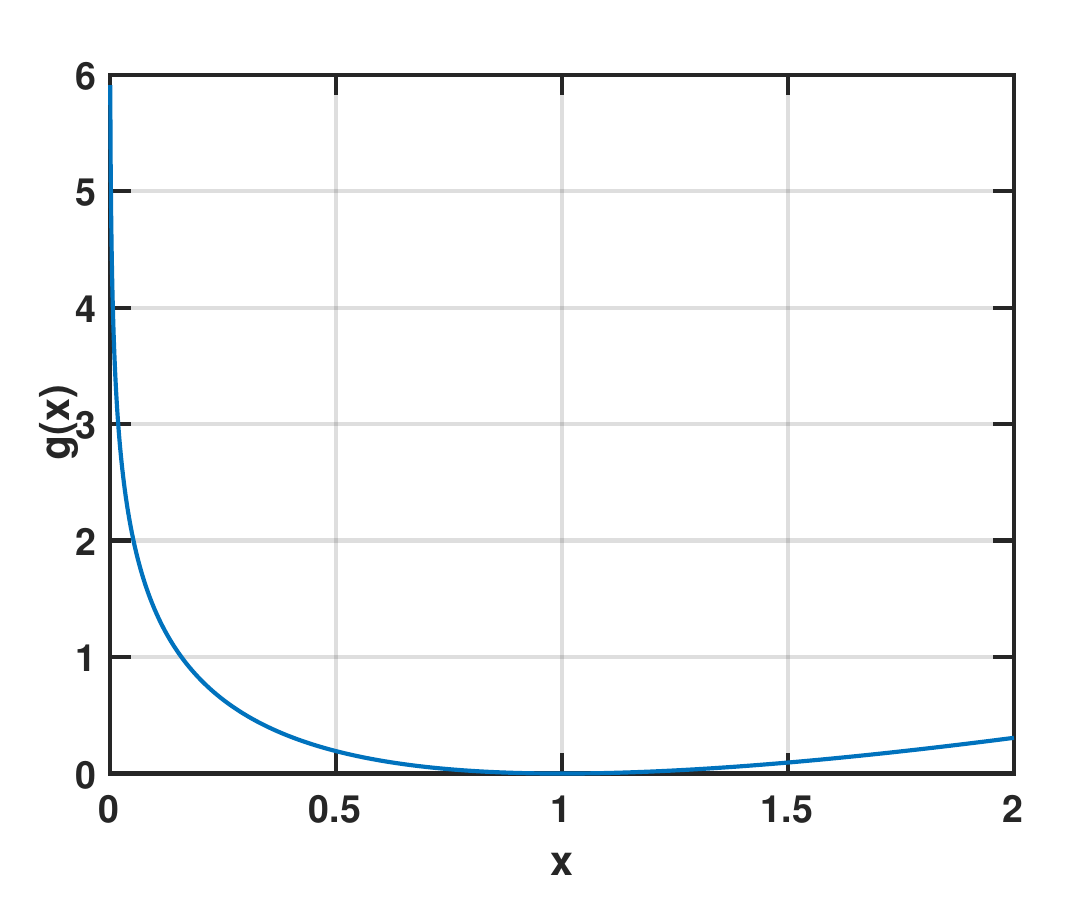}
    \end{center}
    \caption{Plot of $g(x)=x-1-\log(x)$ in the expression of EDD \eqref{eq:EDD_min} with $\delta_\infty^*$; note that $g(x)\geq 0$ which ensures \eqref{eq:EDD_min} is meaningful; $g(x)\rightarrow \infty$ when $x\rightarrow 0$.
    }
    \label{fig:my_label}
\end{figure}

Now, using \eqref{eq:EDD_min} we can further find the optimized window size $w^*$ (for a finite $\gamma$) such that the EDD in  \eqref{eq:EDD_min} is minimized. After taking the derivative with respect to $w$, note that $1/w^2$ is in the order of $o(1/w)$ and thus can be ignored when $w$ is large. Thus we omit all the higher-order terms of $w$ in $D$ and obtain the following result:
\begin{lemma}[Optimal window size and drift parameter]\label{lemma:2}
	For each ARL level $\gamma$, the optimal window size $w^*$ which minimizes the Expected Detection Delay (EDD) in \eqref{eq:final_EDD} is given by:
	\begin{equation}
        \label{eq:optimal_w}
		w^* 
		=\frac{\sqrt{2 \left(1/\Delta-1\right)\left(\sum_{i=1}^m \frac{\lambda_i }{\sigma^2+\lambda_i}B_i\right)}}{mg\left(\Delta\right)}\cdot \sqrt{\log \gamma}.
	\end{equation}
	where $\Delta := \frac 1m \sum_{i=1}^m \frac{\lambda_i}{\sigma^2 +\lambda_i} \in (0, 1)$.
	%
	Substituting the optimal $w^*$ and $\delta_\infty^*$ back into (\ref{eq:d_exp}) gives the optimal value for drift parameter:
	\begin{equation}
		\label{optima_d}
		d^*
		=\frac{m\log(m/D^*)}{m/D^*-1},
	\end{equation}
	where
	\textcolor{black}{$D^*= \sum_{i=1}^m \frac{\lambda_i}{\lambda_i+\sigma^2} \big(1-B_i/w^*+3B_i^2/{w^*}^2\big)$.}
\end{lemma}
\textcolor{black}{The implication of the lemma above is that the optimal window length $w^*$ grows in the order of $\sqrt{\log\gamma}$. Clearly, $w^*$ in \eqref{eq:optimal_w} is positive; in particular, the denominator in \eqref{eq:optimal_w} is non-negative, since $g(x)>0$ for $x\in (0, 1)$ and $\Delta \in (0, 1)$. Moreover, it can be shown that when $w^*$ goes to infinity (since the optimal  $w^*=\sqrt{\log\gamma}$, $\gamma \rightarrow \infty$ in our asymptotic analysis), then $D^*\rightarrow m\Delta$, and $d^*$ in \eqref{optima_d} tends to $\frac{-m\log \Delta}{1/\Delta-1}$.} Figure \ref{fig:winsize} shows numerical examples of choosing the optimized window size $w$. 

\subsection{First-order optimality}

Finally, we show the asymptotic property of Spectral-CUSUM, the EDD ratio relative to the exact CUSUM, under the same ARL constraint; since the exact CUSUM is shown to be optimal and achieves the smallest EDD under a constant ARL constraint as shown in the classic results \cite{lorden1971procedures} and \cite{moustakides1986optimal}. \textcolor{black}{Note that for the ``equalizer'' argument to work, thus we can obtain  ARL and EDD expressions for Spectral CUSUM, we only need the drift term $d$ to be related to $\delta_\infty$ through \eqref{eq:d_exp}. In the following Theorem, we show the resulted EDD under optimal window $w^*$ matches the order of EDD of the Exact CUSUM.
}
%
\begin{theorem}[Asymptotic optimality of Spectral-CUSUM]\label{thm:main}
Given an ARL lower bound $\gamma$, the ratio between the EDD of the Spectral-CUSUM and the EDD of the exact CUSUM satisfies:
	
	\[
	\textcolor{black}{\frac{\mathbb{E}_0[\mathcal{T}_C]}{\mathbb{E}_0[T_C]} = \mathcal O(1)+ 
    \frac{\mathcal I_0 w}{2\log \gamma}\left(1+\mathcal O\left(\frac{1}{\log\gamma}\right)\right).}
	\]
When $\gamma \rightarrow \infty$, for fixed $m$, $n$, $\sigma^2$ bounded away from zero, and $w =\sqrt{\log \gamma}$, this ratio $\mathbb{E}_0[\mathcal{T}_C]/\mathbb{E}_0[T_C]$ tends to an absolute constant, and then the Spectral-CUSUM is asymptotically first-order optimum. 
\end{theorem}
\section{Efficient Computation by  Subspace Tracking}
\label{sec:subspace-track}

In this section, we present an efficient algorithm to keep track of the underlying community structure even if we do not assume the form of the graph structure. This can be treated as a complementary approach to the proposed Spectral-CUSUM method when we have low confidence or lack of pre-change observations. 
Inspired by the GROUSE algorithm \cite{balzano2010online}, which implements stochastic gradient descent on the Grassmann manifold to update subspaces at each time slot, we design a subspace tracking algorithm to update the estimated subspace, denoted as $\hatQ_t$, each time a new graph observation $G_t$ arrives. Here we consider graph observation instead of vector observation $v_t$ on nodes. Similar to \eqref{eq:emp_est}, the graph observation at time $t$ can be represented using vector observation as $G_t = v_t v_t^\top$. 

If we treat each observation as a static graph and apply the spectral clustering method in \cite{von2007tutorial} at each time, we would get a sequence of $\hatQ_t$ independently, which is quite time-consuming and not guaranteed to converge if considering noisy cases. Instead, we can perform an updating procedure each time on $Q\in \mathbb{R}^{n \times m}$ on a Grassmann manifold. The Grassmannian denoted as $\text{Gr}(m,n)$ is a space that contains all $m$-dimensional linear subspaces of the $n$-dimensional vector space. As a compact Grassmann manifold, its geodesics can be computed as indicated in \cite{edelman1998geometry}. Our target matrix $Q$  can then be represented as a point in the Grassmann manifold. Thus the optimization problem becomes finding optimized $\hatQ$ such that:
\begin{equation*}\label{eq:1}
	\hatQ=\arg\min_{Q\in \mathbb{R}^{n \times m}}\sum_{t}\text{tr}(Q^\top G_t Q), ~s.t. ~Q ^\top Q=I.
\end{equation*}
We consider the problem for every time slot and define function $f_t(Q)=\mbox{tr}(Q^\top G_tQ)$ and the derivative of $f_t$ with respect to $Q$ is \cite{petersen2008matrix}:
\begin{equation}\label{eq:3}
	\frac{df_t}{dQ}=\frac{d(\mbox{tr}(Q^\top G_t Q))}{dQ}
	=(G_t+G_t^\top)Q.
\end{equation}
Then we use Equation (2.70) in \cite{edelman1998geometry} to get the gradient of function $f_t(Q)$ on Grassmann manifold from \eqref{eq:3}:
\begin{equation*}\label{eq:4}
	\nabla f_t = (I-QQ^\top)\frac{df_t}{dQ} = (I-QQ^\top)(G_t+G_t^\top)Q.
\end{equation*}

Gradient descent algorithm along a Grassmann manifold is  given by equation (2.65) in \cite{edelman1998geometry}, proving that it is a function of the singular values and vectors of $\nabla f_t$, so suppose we have got the reduced Singular Value Decomposition (rSVD) of $-\nabla f_t=U\Sigma V^\top$ where only the top-$k$  eigenvalues and eigenvectors are kept so that the computational cost is much reduced, we can write the updating function with a step size $\eta$ as:
\begin{equation}\label{eq:5}
	Q(\eta)=\begin{pmatrix}
		QV & U
	\end{pmatrix}\begin{pmatrix}\cos~\Sigma\eta\\\sin~\Sigma\eta\end{pmatrix}V^\top.
\end{equation}
Here we update $Q$ with a step size $\eta$ to get closer to the local minimum on the Grassmann manifold. The complete algorithm is shown in Algorithm \ref{GD}. 
\begin{algorithm}[H]
	\caption{Subspace Tracking for Spectral-CUSUM}\label{GD}
	\begin{algorithmic}[1]
		\Require{ Weighted adjacency matrix at time $t$ denoted as $G_t$, the total number of iterations $T$, number of communities $m$, a set of step sizes $\eta_t$.}
		\Ensure{Subspace Representation $Q$ of the graph. }
		\State Initialize $Q$ randomly,  introduce $y$ as a $T$-dimensional vector with all entries equal to 0.
		\For{$t=1,...,T$}
		\State Observe current adjacency matrix $G_t$.
		\State Compute $\nabla f_t= (I-QQ^\top)(G_t+G_t^\top)Q$.
		\State Compute SVD of $-\nabla f_t=U\Sigma V^\top$.
		\State Update $Q$ using \eqref{eq:5} with step size $\eta_t$.
		\State Update $y(t) = \text{tr}(Q^\top G_t Q)$
		\EndFor
		\State Run CUSUM detection procedure on $y$ and get detection statistic $S$.
		\State \Return {$S$}
		
	\end{algorithmic}
\end{algorithm}

{\it Choice of Step Size.} For our problem, constant step size and decreasing step size can be efficient. The constant step is slower initially, but it is more stable to detect changes in community structures. In our numerical experiments, we take constant step size $\eta= 0.01$.

{\it Complexity.} In practice, we notice that in this algorithm, we have to perform SVD for each iteration, which is quite time-consuming. To accelerate the whole algorithm, we use the incremental SVD algorithm \cite{sarwar2002incremental} in each iteration.


\section{Numerical Experiments}
\label{sec:exp}
In this section, we compare our method with the state-of-the-art and discuss their numerical results on both synthetic and real data sets. Since a graph can be considered as a discrete approximation to a manifold \cite{bishop2011geometry}, we also show that our model can achieve promising performance on dynamic manifold data.

\subsection{Methods for comparison} 
In our experiments, we compare our method with four other baseline approaches, including (1) Generalized Likelihood Ratio procedure based on vectorized data (\texttt{Vectorized GLR}) (2) Hotelling's $T$-squared CUSUM (\texttt{Hotelling}); (3) Single eigenvector procedure (\texttt{SC,$m=1$}); (4) Rank-$m$ Subspace Tracking of Section \ref{sec:subspace-track} (\texttt{SGD}); and (5) Exact CUSUM as a sanity check. 
The detailed explanation of these baseline methods is as follows:
 
(1) {\it Generalized Likelihood Ratio (GLR) procedure based on vectorized data.} This baseline method completely ignores the topology properties of the adjacency graph $G_t$ and vectorizes it as $g_t = \text{vec}(G_t)$ such that the previous hypothesis test (\ref{eq:model}) becomes: 
\begin{equation}
	\label{eq:vec_model}
	\begin{array}{ll}
		H_0:  g_t \distas{\text{iid}} \mathcal{N}(0,\sigma^2I_{n^2})& t=1,2,\dots,\tau\\
		H_1:  g_t \distas{\text{iid}} \mathcal{N}(h,\sigma^2I_{n^2})~& t = \tau+1,\tau+2,\dots\\
	\end{array}
\end{equation}
where $h=\text{vec}(AA^\top)$. Then classical result from \cite{lorden1971procedures,lai-ieeetit-1998} gives the optimal stopping criteria for GLR procedure:
\begin{equation}
	\label{eq:vec_GLR}
	T_{\text{GLR}} = \inf\bigg\{t:\max_{t-w<k<t}\frac{(\sum_{i=k+1}^t \left\Vert g_i \right\Vert)^2}{t-k}  > b\bigg\}.
\end{equation}

(2) {\it Hotelling's $T$-squared CUSUM.}  We also compare our method with Hotelling's $T$-squared statistic, which was introduced in \cite{hotelling1992generalization}. The way to implement this is to calculate a pre-change sample mean $\hat{\mu}_0$ and sample covariance $\widehat{\Sigma}_0$, then keep track of the current sample mean $\bar\mu_{t-w,t}$ via a sliding window. We also need to pre-define the drift parameter $d^{H_2}$, which can be obtained from historical data. Thus the detection statistic is given by:
\[
S^{H_2}_t = (S^{H_2}_{t-1})^+ + (\bar\mu_{t-w,t}-\hat{\mu}_0)^\top \widehat{\Sigma}_0^{-1}   (\bar\mu_{t-w,t}-\hat{\mu}_0) - d^{H_2} .
\]
Then the optimal stopping time for Hotelling's statistics is given by:
\begin{equation}
	\label{eq:hotelling}
	T^{H_2} = \inf\bigg\{t>0:S^{H_2}_t \geq b\bigg\}.
\end{equation}

(3) {\it Single eigenvector procedure of Spectral-CUSUM.} This baseline method is almost the same as the proposed Spectral-CUSUM with the only difference of fixing community size $m=1$. We use it to illustrate the importance of finding the optimal community size. 

(4) {\it Rank-m Subspace Tracking.} This is proposed in Section \ref{sec:subspace-track} as a numerical alternative method for our Spectral CUSUM. So we list it as another baseline. 


For synthetic data, we study the relationship between ARL and EDD. For real data, we study EDD only since we cannot implement simulation to decide the average run length if the history data under pre-change distribution is not sufficient.


\subsection{Synthetic experiments}

We devise three synthetic experiments to study the performance of our method on the emergence and the switching membership problems.




\vspace{.1in}
\emph{Synthetic data 1} is used for evaluating our method of detecting emerging communities. In this experiment, we set $\sigma=5$ and assume a dynamic graph that contains 50 nodes without any community structure in the beginning. 
When the change occurs at time $t=$, three communities are formed, containing 10, 10, and 15 nodes, respectively. In other words, we have $m=0$ before $t$ and $m^*=3$ after $t$. 
The optimal window size can be found according to \eqref{eq:optimal_w}, where $w=5$.

\begin{figure}[ht!]
	\begin{center}
		\includegraphics[width = 0.35\textwidth]{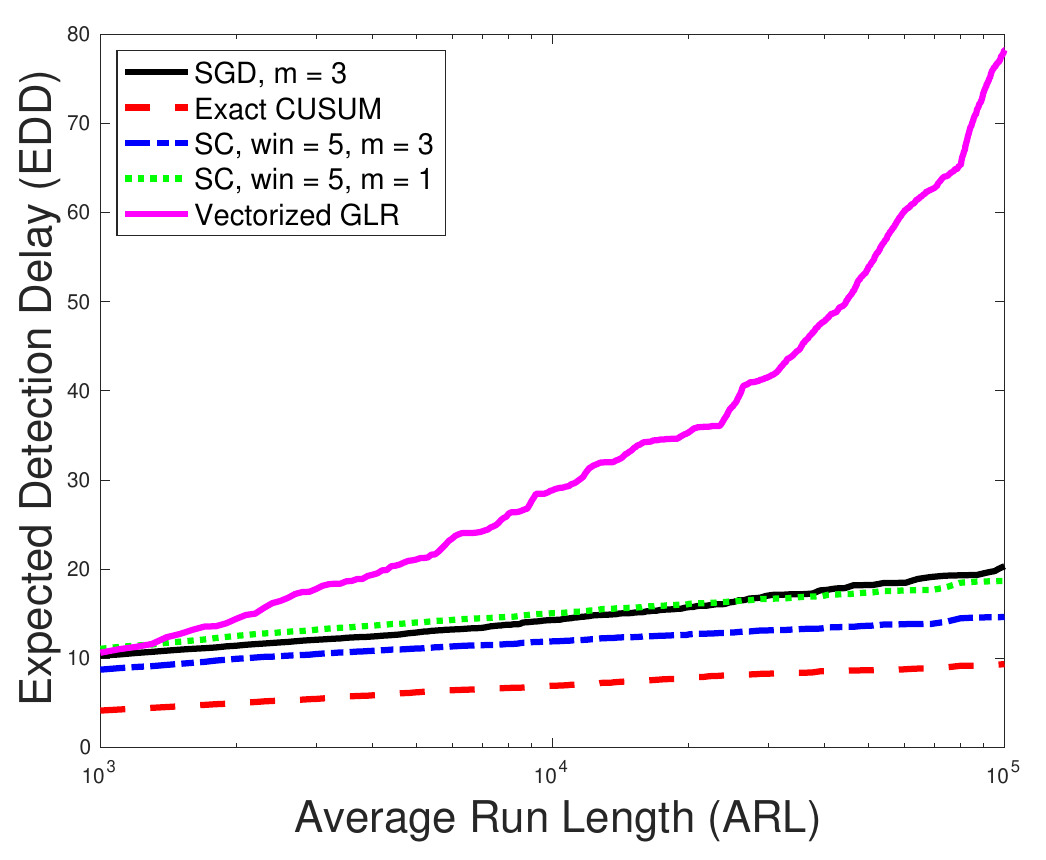}
	\end{center}
	\caption{Comparison of the single eigenvector procedure and
		Exact-CUSUM procedure for emergence problem. Community size fixed and $\sigma$ is set to be 5.}
	\label{fig:ARLEDD}
	
\end{figure}

\begin{figure}[ht!]
	\begin{center}
		\subfigure[Smaller emerging communities]{
			\includegraphics[width = 0.35\textwidth]{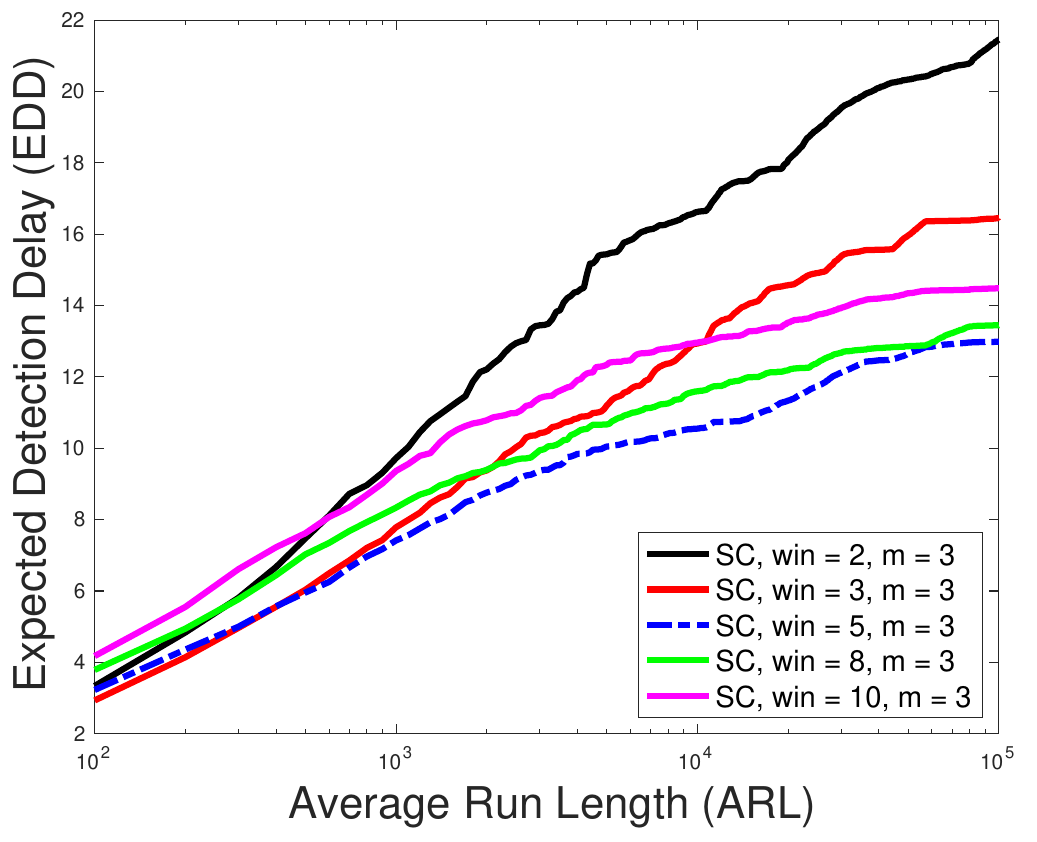}}
		\subfigure[Larger emerging communities]{
			\includegraphics[width = 0.35\textwidth]{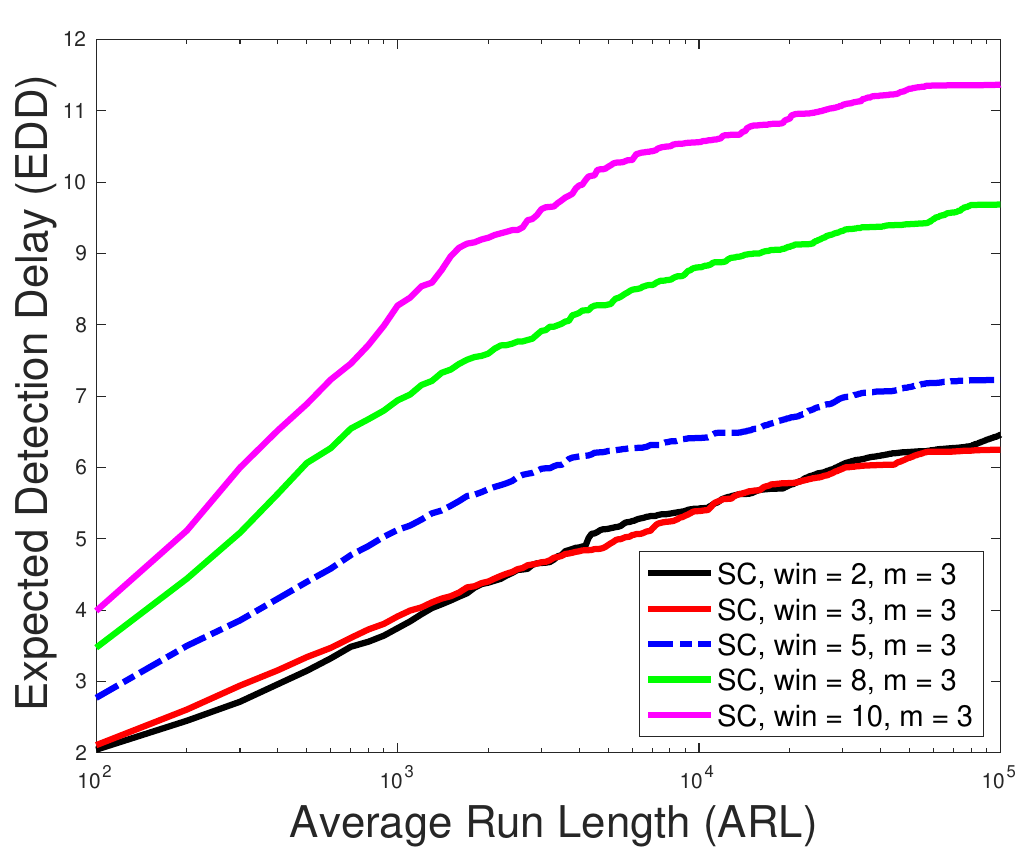}}
	\end{center}
	\caption{ARL vs. EDD plot for optimized window size for Spectral-CUSUM with smaller (a) and larger (b) emerging communities.}
	\label{fig:winsize}
	
\end{figure}

\begin{figure}[ht!]
	\begin{center}
		\includegraphics[width = 0.35\textwidth]{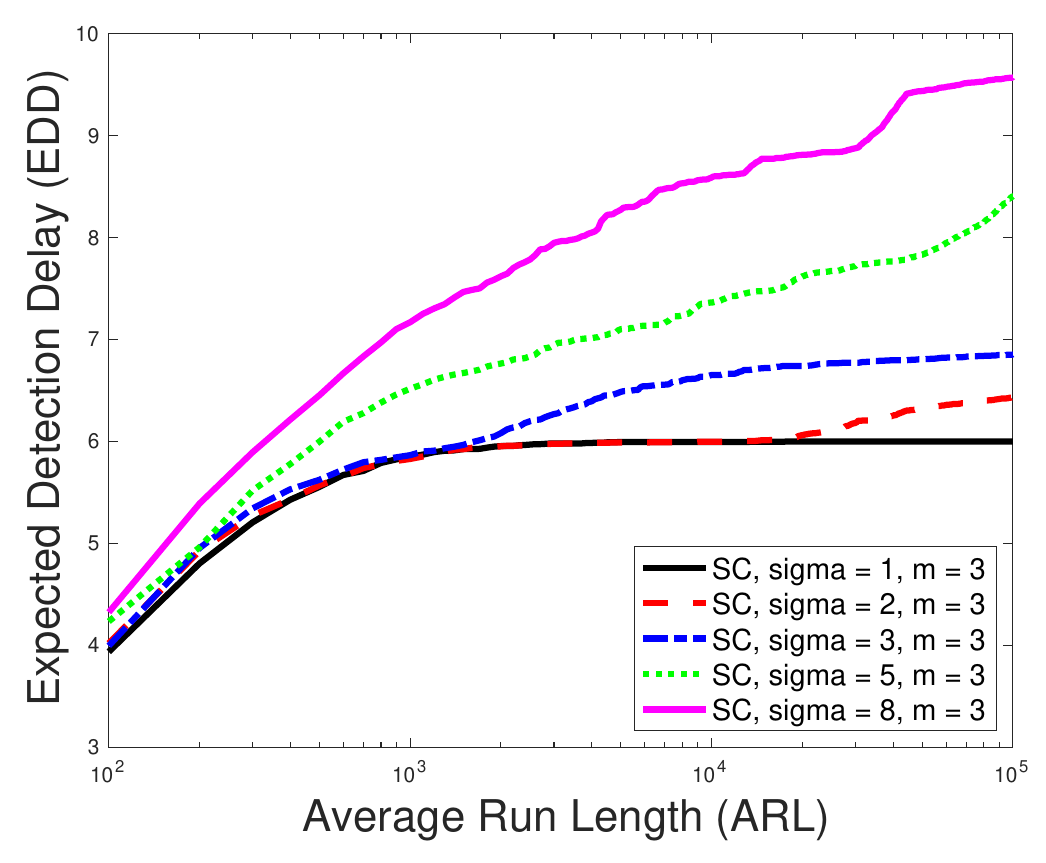}
	\end{center}
	\caption{Comparison of the different noise levels for emergence problem. Community size is fixed, and values of $\sigma$ are set to be 1,2,3,5,8, respectively.}
	\label{fig:sigma_comp}
	
\end{figure}

The experimental result shows that our method significantly outperforms other baselines. 
Figure~\ref{fig:ARLEDD} presents the ARL of all the methods on the synthetic data 1 with different choices of community density and window size. 
As can be observed, our method (blue line) works much better (the lower, the better) than the baselines and is closest to the sanity check method Exact-CUSUM. A comparison of how community size would affect the choice of optimal window size is also shown in the figure~\ref{fig:winsize}$(b)$ where we increase the size of each emerging community by 5, referred to as larger emerging communities. It can be inferred that the change becomes easier to detect in larger communities. Thus, the optimal window size decreases correspondingly. In figure~\ref{fig:sigma_comp}, we show how different noise levels could infect the general detection delay, which leads to the conclusion that a larger noise level leads to a longer average detection delay.  Moreover, when ARL is relatively small, a large noise level would not have much effect on the EDD. However, if the ARL is larger, the EDD will be increased significantly for higher noise levels.

%

	

\begin{figure}[ht!]
	\begin{center}
		\includegraphics[width = 0.35\textwidth]{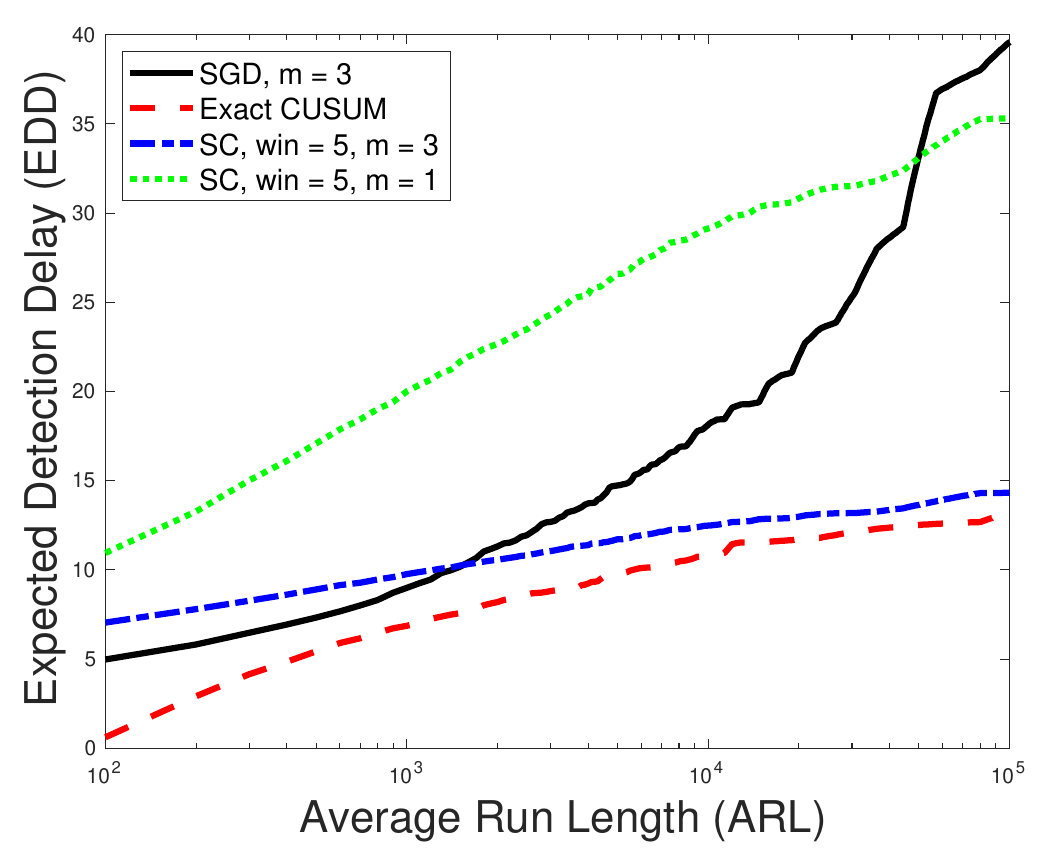}
	\end{center}
	\caption{Comparison of the single eigenvector procedure and Exact-CUSUM procedure for switching membership problem. Community Size fixed to be 3. The total number of nodes in the network is set to 50. Before the change, each community has a size of 10, and after the change, each community's size increases to 15 by absorbing nodes that do not belong to any community before the change happens.}
	\label{fig:switch_ARL}
	
\end{figure}

\vspace{.1in}
\emph{Synthetic data 2} is used for evaluating our method of detecting switching membership. For the design of the switching membership experiment, the number of nodes is 50, and the size of the three communities remains the same: 10, 10, and 15 nodes. However, the membership for all three communities will change. A comparison between SGD, Exact-CUSUM, and Spectral-CUSUM is shown in Figure \ref{fig:switch_ARL}. We can see that the single eigenvector procedure ($m=1$) performs much worse than Spectral-CUSUM using the correct potential community size ($m=3$). The performance gain of spectral-CUSUM in the switching membership setup is larger than that in the emerging subspace problems compared with the Exact-CUSUM procedure. Moreover, the EDD of SGD is growing quite fast with the increase of ARL, which indicates that it is noise-sensitive. 

Here we do not compare with the vectorized GLR procedure since it does not work in the switching membership problem. The reason is that the switching membership change does not lead to the increase or decrease of average edge weights; thus, the change cannot be detected without considering the graph topology.

\begin{figure}[!ht]
	\begin{center}              
		\includegraphics[width = 0.35\textwidth]{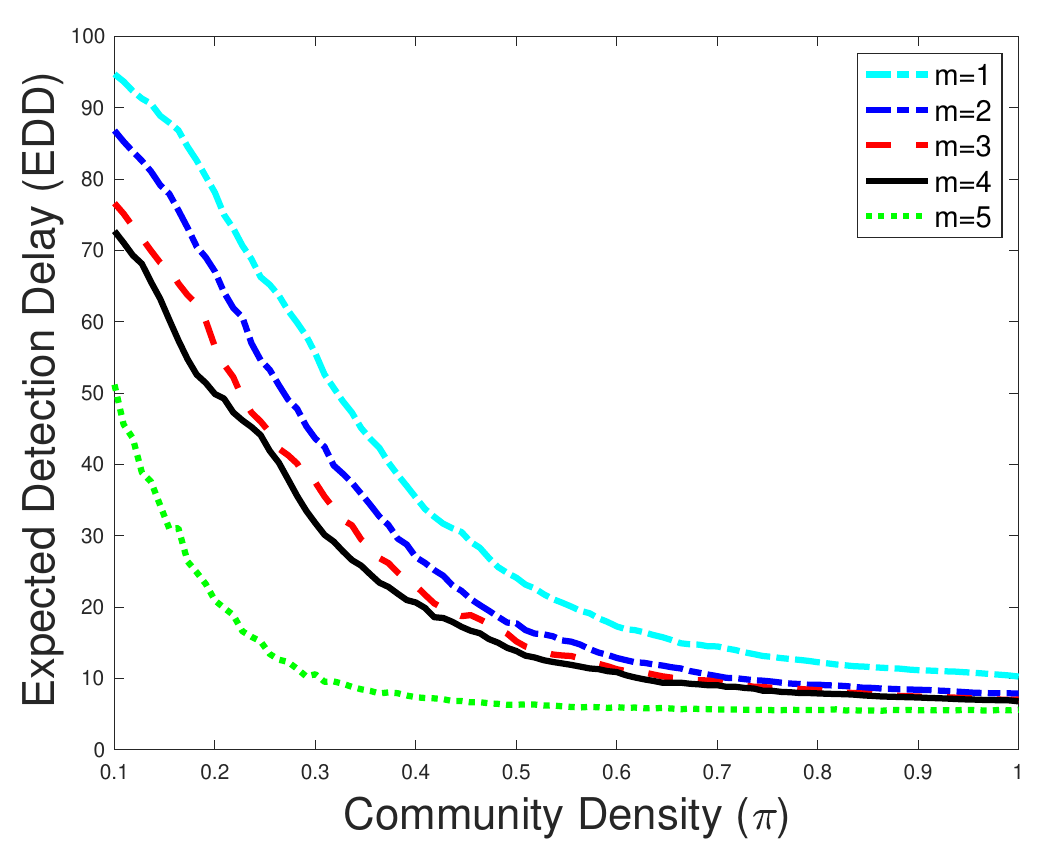}
	\end{center}
	\caption{Relationship between the community density and EDD under different settings of the number of communities  $m=\{1,2\dots,5\}$. for emergence community problem. Smaller $m$ leads to higher EDD, and a larger community density $\pi$ after the change would have also decreased the EDD.}
	\label{fig:pVSEDD}
\end{figure}

\vspace{.1in}
\emph{Synthetic data 3} is used to study the impact of network topology. We also explore the influence of community density by assuming that for community $i$, the probability of forming an edge within the community is $p_i$. Then by changing the value of $p_i$ under an emergence community scenario, the density of the community can change. From Figure \ref{fig:pVSEDD} we see that a denser community structure after the change leads to a more accurate detection. However, we can also observe from the drastically decreasing curve that the detection can be performed more accurately if the community density $p_i \geq 0.4$. This provides an approach to discovering a community through the Spectral CUSUM method. In addition, we can see that the detection power increases with a larger  $m$ and the increase in community density. This result is consistent with our intuition since we can always collect more useful information when $m\leq k$, and the denser the community is, the more signal that reveals the structure would be much stronger, leading to quicker change detection.


\subsection{Detecting changes of manifold structure}

\begin{figure}[!ht]
	\begin{center}
		\subfigure[before change]{              
			\includegraphics[width = 0.22\textwidth]{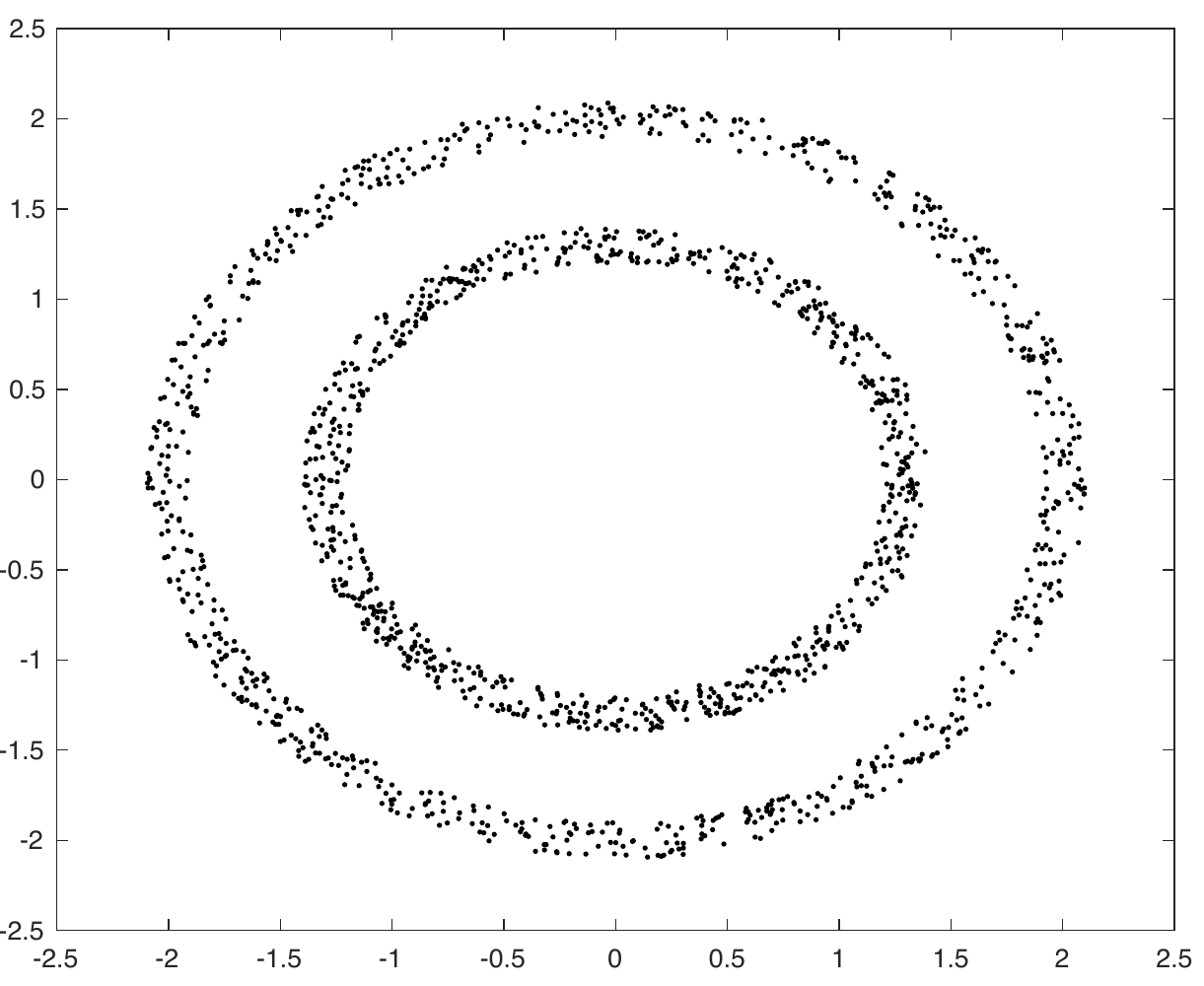}}
		\subfigure[after change]{
			\includegraphics[width = 0.22\textwidth]{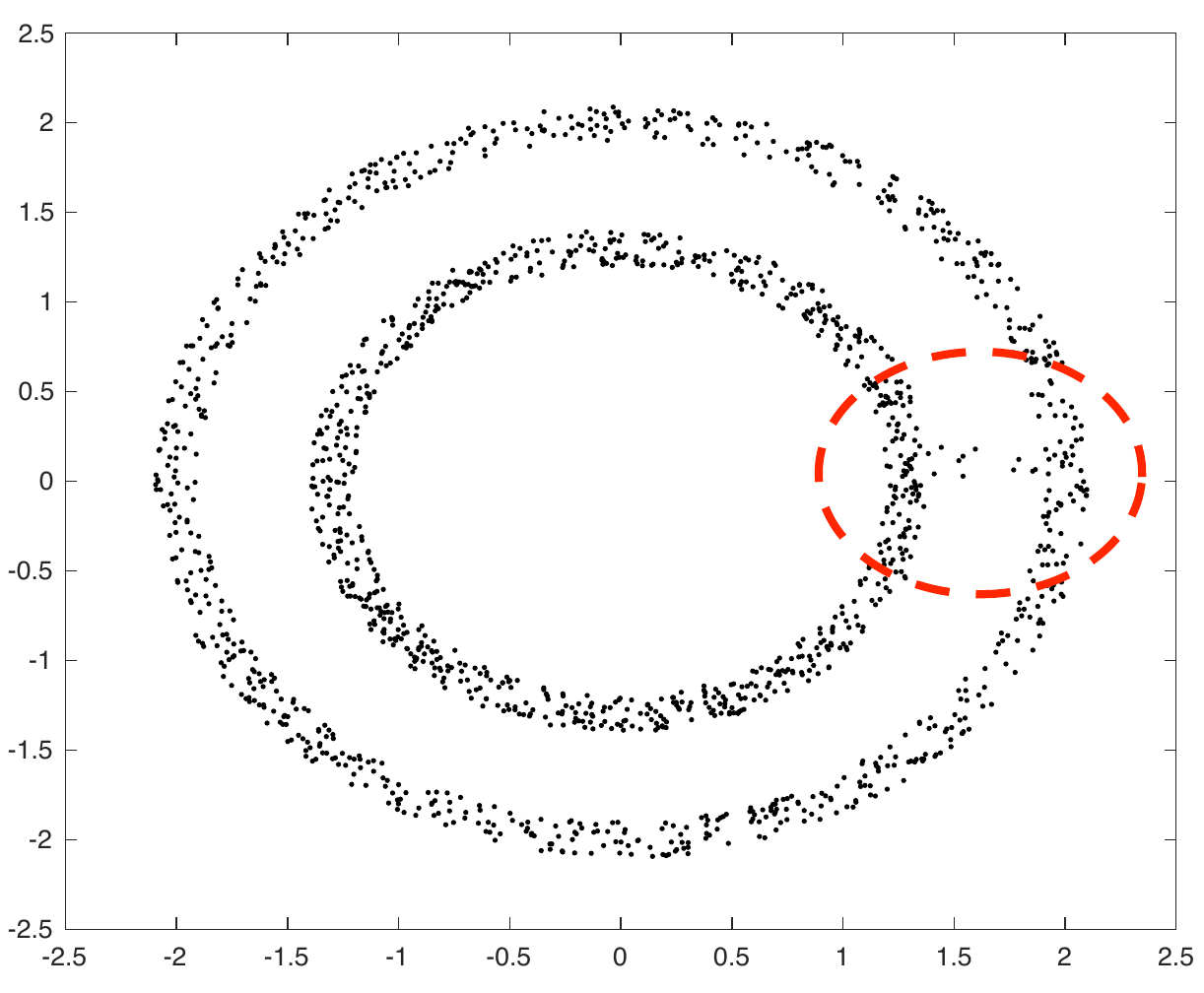}}
		\subfigure[ARL VS EDD plot]{
			\includegraphics[width = 0.35\textwidth]{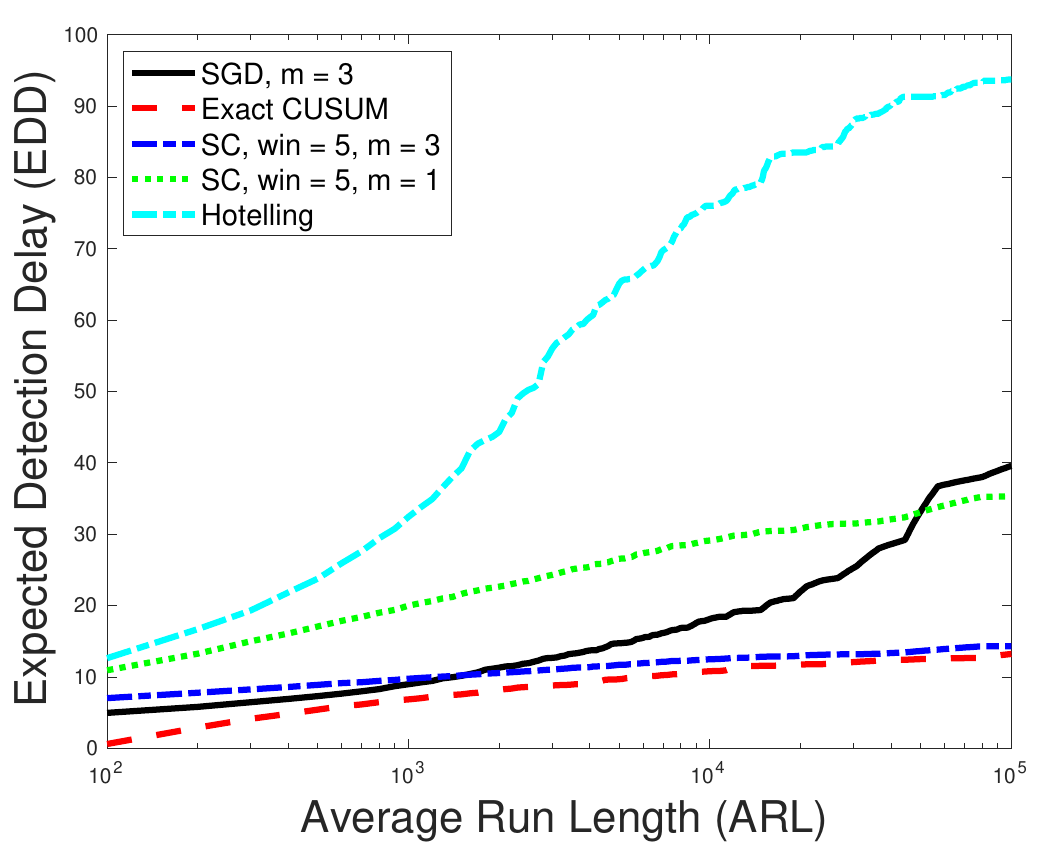}}
	\end{center}
	\caption{A case that two separate rings have some overlapping nodes by a bridge after the change: (a) Before the change, there are two rings structures that form two communities; (b) After the change, two communities remain the same, but some of their members join another community, leading to the overlapping of two communities. (c) The ARL VS EDD results show that Spectral-CUSUM when $m=4$ outperforms other procedures under larger ARL, but SGD is better for smaller ARL.}
	\label{fig:manifold_illu}
\end{figure}

This section shows that our detection algorithm also works for detecting structural changes in the manifold. We first consider a one-dimensional manifold on two-dimensional rings. Before the change happens at time $t=200$, the two rings are separately located, while after the change, there forms a bridge between the two rings, as shown in Figure \ref{fig:manifold_illu}. 

We use the Isomap \cite{balasubramanian2002isomap} to calculate the adjacency graph, which is sparse. Given the adjacency matrix, we set the threshold to eliminate all the edges whose distances run below the threshold. Thus, inside one community, any pair of members can be reached via a few steps. Using the Spectral-CUSUM, we show that, in this case, we can detect the change quite quickly in figure~\ref{fig:manifold_illu}(c).

\begin{figure}[!ht]
	\label{fig:swissroll}
	\begin{center}
		\subfigure[before change]{
			\includegraphics[width=0.23\textwidth]{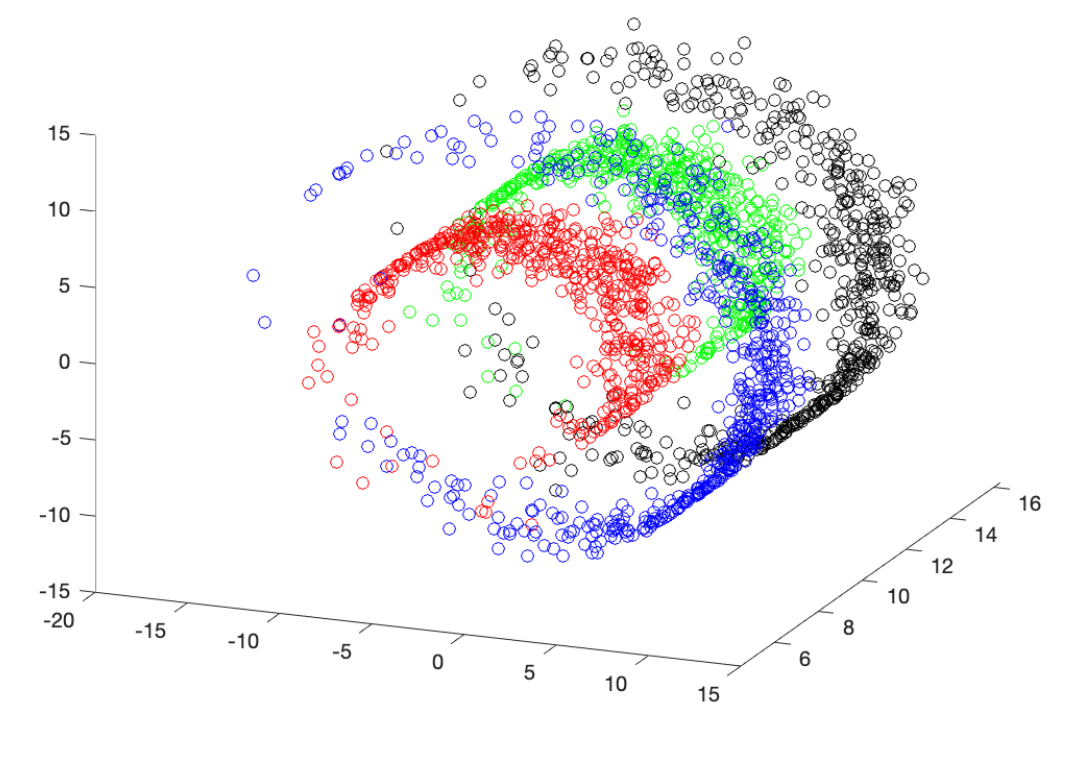}}
		\subfigure[after change]{
			\includegraphics[width=0.23\textwidth]{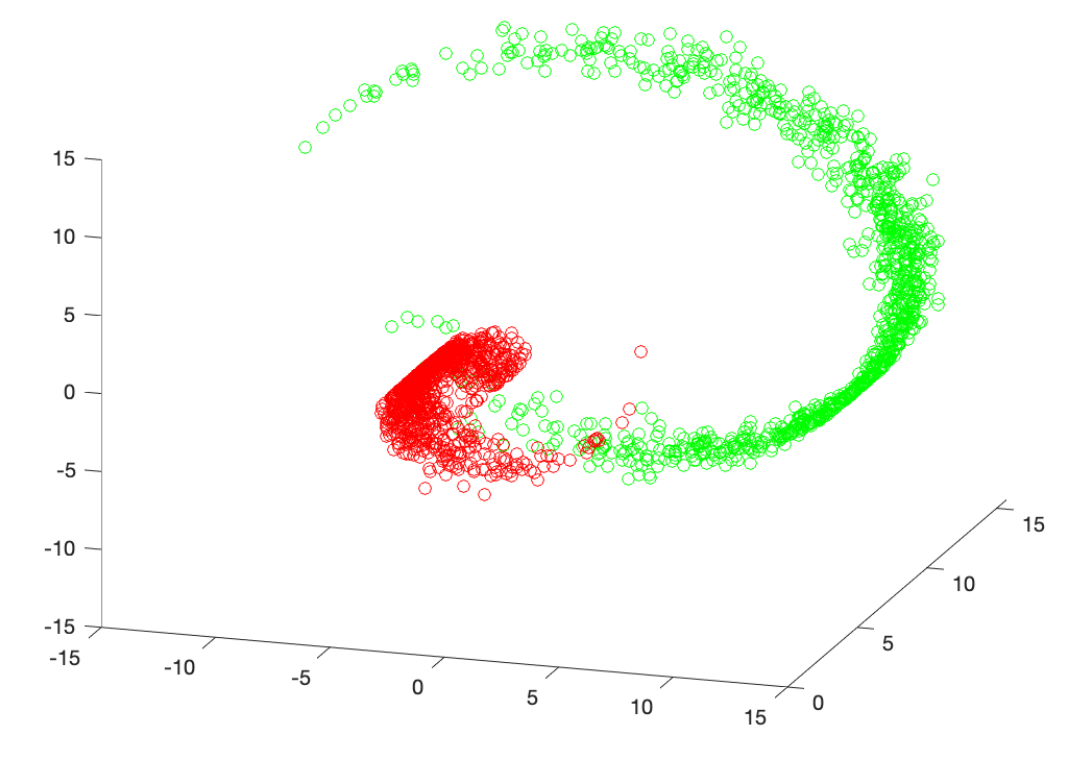}}
		\subfigure[ARL VS EDD plot]{
			\includegraphics[width = 0.35\textwidth]{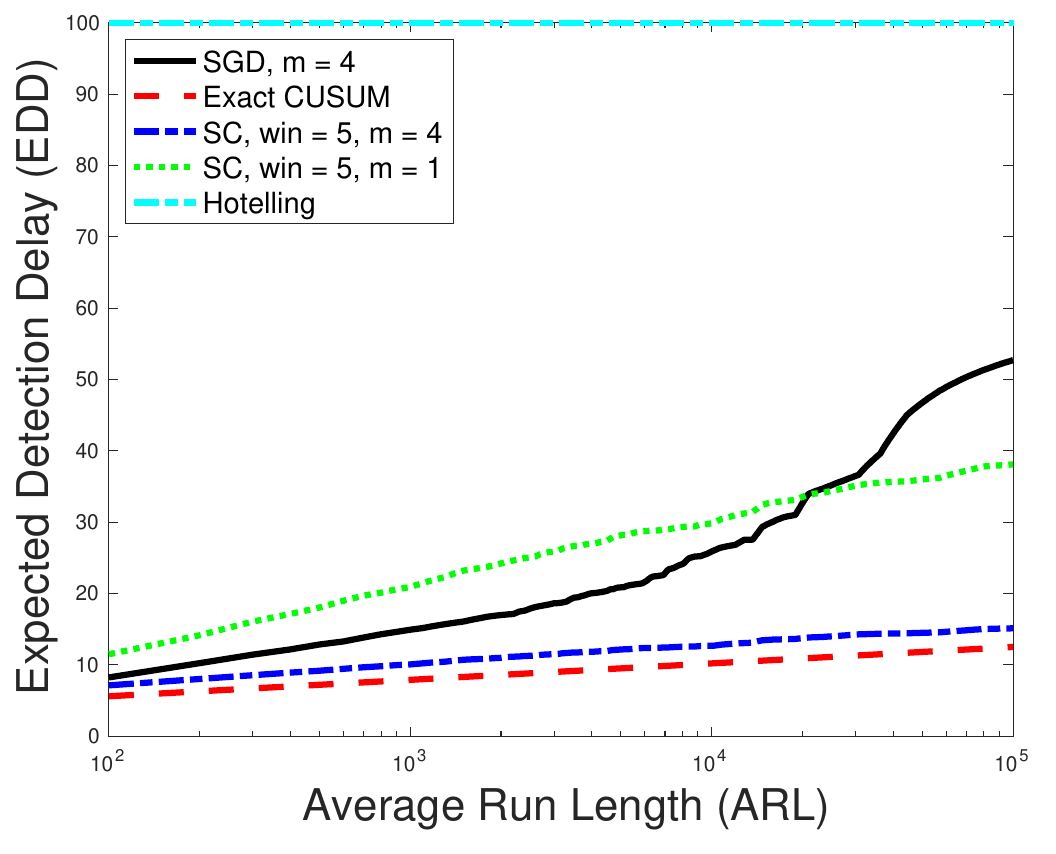}}
	\end{center}
	\caption{A case that four Swiss rolls merge into two communities after the change: (a) Before the change, there are four swiss roll structures which form four different communities; (b) After the change, two communities disappear and merge into the other communities. Note that in this experiment, Hotelling's $T$-squared statistics cannot detect the change at all.}
\end{figure}

Another classical manifold is Swiss Roll data \cite{moolenbeek1981swiss}, which is a two-dimensional structure lying in a three-dimensional space. The experiment is designed such that four Swiss rolls merge into two Swiss rolls after the change, as shown in Figure \ref{fig:swissroll}. We aim to use this to show that our method can deal with various data topological structure changes based on a choice of similarity measure. We see from the result that when $m$ equals its potential community size, the performance is quite close to the optimal Exact-CUSUM. 

\subsection{Yellowstone seismic sensor network data}

\begin{figure}[!ht]
	\begin{center}
		\subfigure[$t=10$]{              
			\includegraphics[width = 0.23\textwidth]{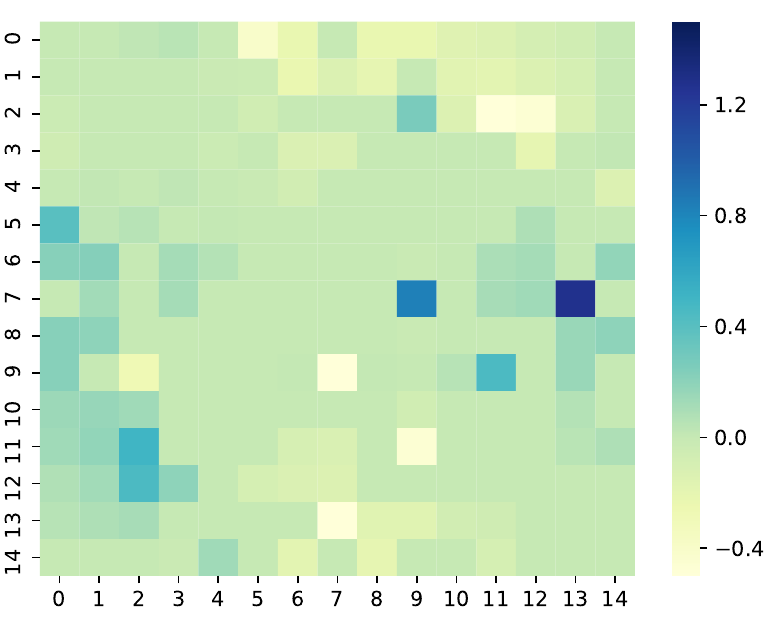}}
		\subfigure[$t=180$]{
			\includegraphics[width = 0.23\textwidth]{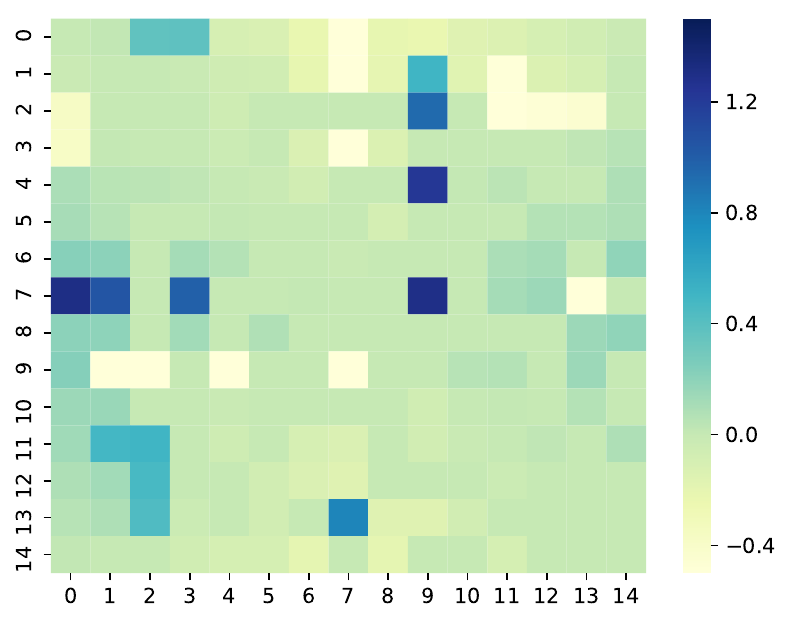}}
		\subfigure[$t=205$]{
			\includegraphics[width = 0.23\textwidth]{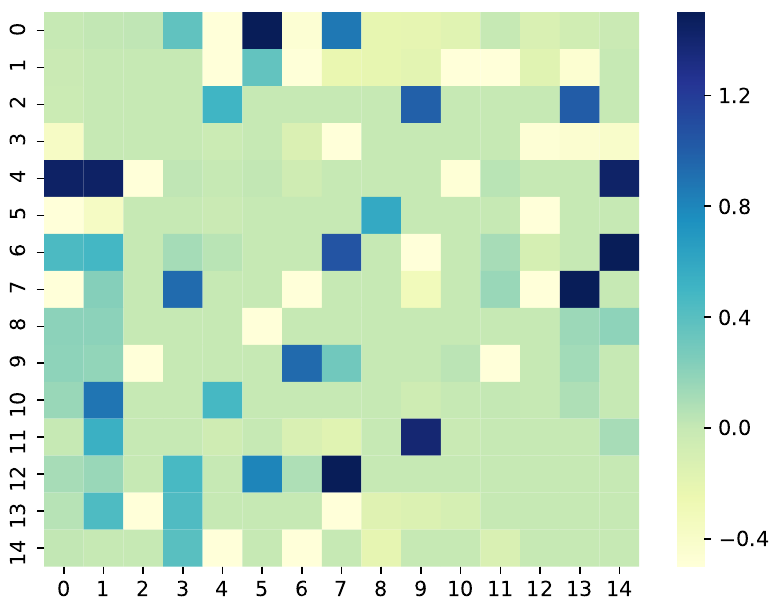}}
		\subfigure[$t=250$]{
			\includegraphics[width = 0.23\textwidth]{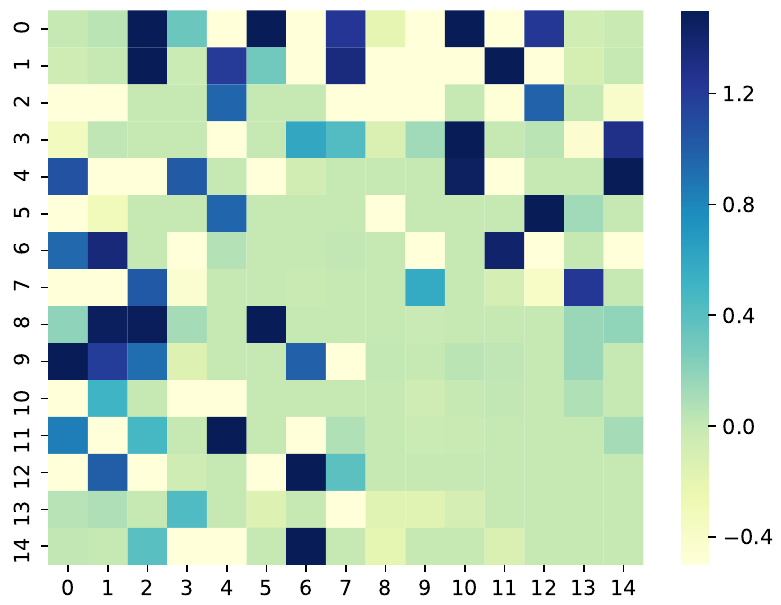}}
	\end{center}
	\caption{Correlation matrix for seismic sensors at different times $t=10,180,205,250$. We can see clearly from the above matrices that there's a big change before and after $t=200$, which is the ground-truth change-point in the seismic network, causing an emergence community scenario.}
	\label{fig:seismic_data}
\end{figure}

We further consider a seismic sensor network data set adapted from \cite{he2018sequential}. The sensors are placed in different locations to measure signals around the Old Faithful Geyser in Yellowstone National Park. The total number of sensors is 19 in this case, but we removed four sensor signals since they failed to work during data collection. We then observe a sequence of sensor signals and translate each one into a dynamic cross-correlation graph. At the very beginning, the cross-correlation between each pair of sensors is low, which means they are not related. The emergence of community happens in the middle of the sequence. Such change will cause the infected sensors to generate similar signals, which lead to a higher correlation magnitude between them. In contrast, an unaffected sensor still generates random noises, thus having a low correlation with other sensors. Consequently, a community containing all the affected sensors emerges after the change happens. We visualize the correlation matrix for better understanding in Figure \ref{fig:seismic_data}. In this case, the true changepoint time is roughly known, corresponding to the geyser eruption time. 

\begin{figure}[!ht]
	\begin{center}   
		\subfigure[Seismic topology visualization]{           
			\includegraphics[width = 0.3\textwidth]{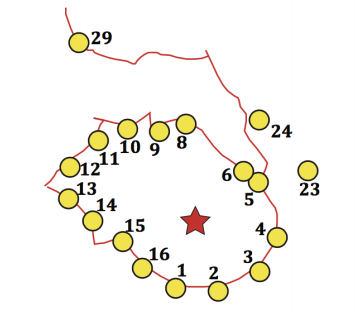}}
		\subfigure[Detection statistics]{
			\includegraphics[width = 0.39\textwidth]{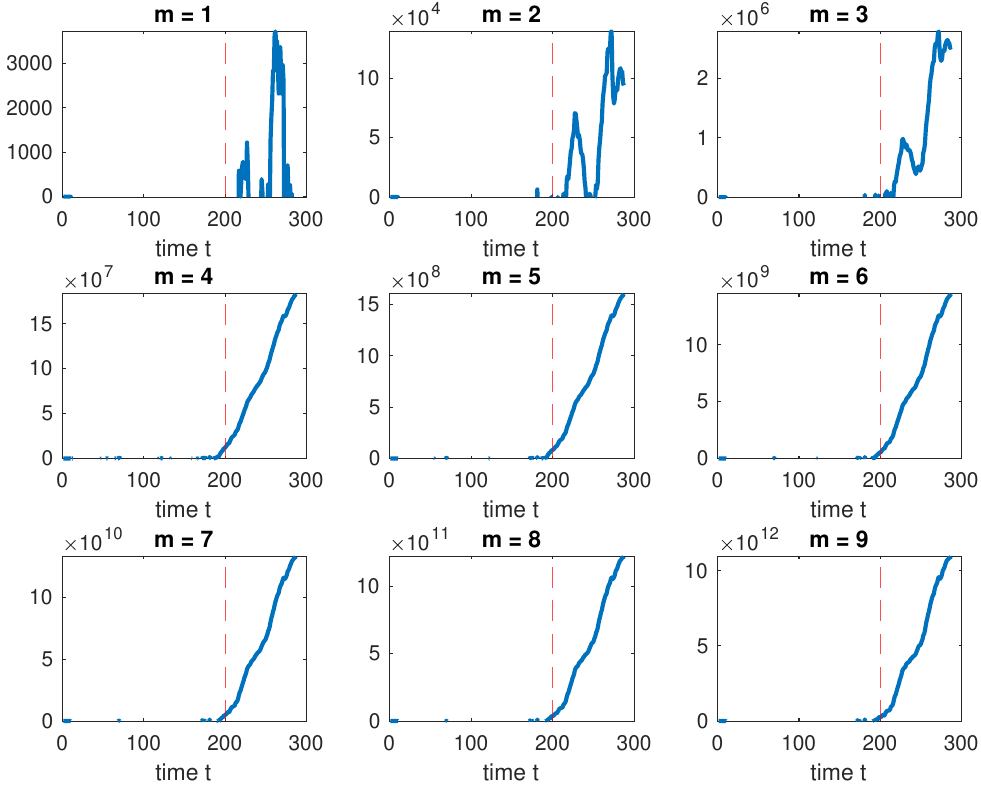}}
		
	\end{center}
	\caption{The topology of seismic sensors and detection statistics for sequential seismic data with an outburst of community structure changes. In (a), the total number of 19 seismic sensors roughly forms a circle, and signals can be observed on each sensor \cite{wu2017anatomy}. Thus the community would be formulated when an earthquake happens, making several sensors correlated. In (b), we applied our Spectral-CUSUM statistics on sensor signals, and we know the ground-truth earthquake happens at time 200. The results demonstrate the ability of the algorithm to detect such changes very quickly.}
	\label{fig:seismic}
	
\end{figure}

We apply our Spectral-CUSUM procedure to the dataset with choices of different potential community sizes $m$. Figure \ref{fig:seismic} (a) shows the sensor locations, and the results of detection statistics using Spectral-CUSUM procedure are shown in Figure \ref{fig:seismic} (b). It can be seen clearly that setting potential community size $m\geq4$ gives better performance on detecting change at around $t=200$, which corresponds to a geyser eruption time -- treated as a true change-point time. The result shows that the event can be detected without false alarms when $m\geq 4$, which reveals that this seismic network's true underlying community size is about 4. The detection delay is shown in Table \ref{table1}. This shows that our detection statistics can be useful when the underlying graph structure is unknown.

\begin{table*}[ht]\small\setlength\tabcolsep{10pt}
	\renewcommand\arraystretch{1.2}
	\centering
	\caption{\small Detection delays under different $\#$ communities $(m)$ for seismic data.}
	\begin{tabular}{cccccccccc}
		\specialrule{.08em}{0em}{0em}
		$\#$ communities $(m)$  & 1  & 2  & 3  & 4 & 5 & 6 & 7 & 8 & 9 \\ 
		\hline
		Detection Delay & 23 & 20 & 17 & 4 & 3 & 3 & 2 & 2 & 2 \\ 
		\specialrule{.08em}{0em}{0em}
	\end{tabular}
	\label{table1}
\end{table*}

\begin{figure}[!ht]
	\begin{center}              
		\includegraphics[width = 0.35\textwidth]{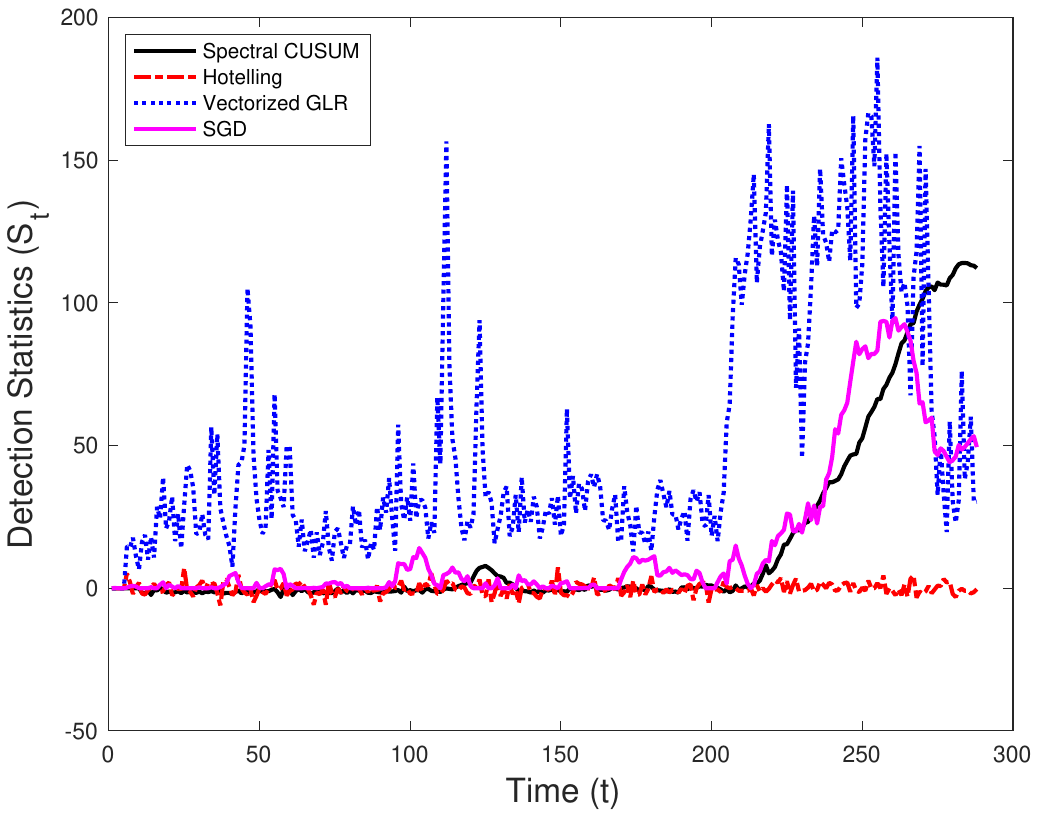}
	\end{center}
	\caption{Comparison of our Spectral-CUSUM procedure with other baseline methods, including Hotelling's statistics, SGD, and vectorized GLR.}
	\label{fig:real_base}
\end{figure}

In addition, we compare our Spectral-CUSUM with other methods on this real data. We can see from Figure \ref{fig:real_base} that Spectral-CUSUM and SGD have the most obvious abrupt change around the true change-point time 200. However,  vectorized GLR and Hotelling's $T$-squared statistics can not detect such a change.

\section{Conclusion}\label{sec:conclusion}

We present a novel Spectral-CUSUM procedure for detecting underlying community changes through noisy observations. We provide the first-order asymptotic optimality of Spectral-CUSUM under the optimal choice of parameters. We also present an efficient online computation procedure to evaluate the Spectral-CUSUM statistic without remembering all data based on subspace tracking. Experimental results on both synthetic and real-world data demonstrate the superior performance of our procedure compared with alternative methods.


%

\section*{Acknowledgment}

The work is partially supported by National Science Foundation CAREER Award CCF-1650913,  NSF CMMI-2015787, DMS-1938106, and DMS-1830210. 

\ifCLASSOPTIONcaptionsoff
  \newpage
\fi



%
\bibliographystyle{IEEEtran}
\bibliography{refs}


\newpage
\appendix

\begin{proof}[Deriving the bounds for terms in \eqref{u_est_error}] Recall $\rho_i=1/(\sigma^2 + |C_i|)$. 
First, the Frobenius norm for the second term in the covariance matrix of the estimation error:
\begin{align*}
	& \frac{\rho_i/\sigma^2}{(\rho_i-1/\sigma^2)^2} \|I- UU^\top\|_F \\
	= &~ \frac{\sigma^2/\rho_i}{\left(\sigma^2 - 1/\rho_i\right)^2} \|I- UU^\top\|_F\\
	= &~\frac{\sigma^2(\sigma^2+|C_i|)}{|C_i|^2} \|I- UU^\top\|_F \\
	= &~ \frac{\sigma^2}{|C_i|}\left(\frac{\sigma^2}{|C_i|}+1\right) \sqrt{n-m},
\end{align*}
where we have used $\|I- UU^\top\|_F^2 =\mbox{tr}[(I-UU^\top)(I-UU^\top)] =\mbox{tr}(I-UU^\top)= n-m$ since $U$ is a $n$-by-$m$ orthonomal matrix and $I-UU^\top$ is a projection matrix. Assume $\sigma^2/|C_i| < \epsilon$, $\epsilon > 0$, we have the last equation above is upper bounded by $\epsilon(1+\epsilon)\sqrt{n-m}$. 

On the other hand, it can be shown that the square Frobenius norm of the first term $\sum_{k=1,k\neq i}^{m}\frac{\rho_i \rho_k}{(\rho_i-\rho_k)^2}u_ku_k^\top$ is given by (and further bounded by)
\begin{equation}
\begin{split}
    &\sum_{k=1,k\neq i}^{m}\frac{\rho_i^2 \rho_k^2}{(\rho_i-\rho_k)^4} \\
    = & \sum_{k\neq i} \frac{\frac{1}{\rho_i^2}\frac{1}{\rho_k^2}}{\left(\frac{1}{\rho_k}-\frac{1}{\rho_i}\right)^4}\\
    = & \sum_{k\neq i} \frac{(\sigma^2 + |C_i|)^2(\sigma^2 + |C_k|)^2}{(|C_k| - |C_i|)^4}\\
    =& \sum_{k\neq i} \frac{\left(\frac{\sigma^2}{|C_i|} + 1\right)^2\left(\frac{\sigma^2}{|C_k|} + 1\right)^2}{\left(
    \sqrt{\frac{|C_k|}{|C_i|}} - 
    \sqrt{\frac{|C_i|}{|C_k|}}
    \right)^4} \\
    \geq & \textcolor{black}{\frac{m-1}{\left(
    1+\theta - \frac{1}{ 1+\theta}
    \right)^4}}, 
    \end{split}
\end{equation} 
where the last inequality is based on Assumption \ref{assumption}(2). 

Thus, to have the first term dominate the second term, we need 
\[
\epsilon(1+\epsilon) <\eta \cdot \sqrt{\frac{m-1}{n-m}}\cdot \textcolor{black}{\frac{1}{\left(
    1+\theta - \frac{1}{ 1+\theta}
    \right)^2}}, 
\]
where $\eta \in (0, 1)$ controls the relative magnitude of the terms, which gives us the condition 
\[
\epsilon <  
\left(\frac 1 4 +
\eta \cdot \sqrt{\frac{{m-1}}{{n-m}}}\cdot \textcolor{black}{\frac{1}{\left(
    1+\theta - \frac{1}{ 1+\theta}
    \right)^2}}
\right)^{1/2} - \frac{1}{2}.
\]
Together with $\frac{\sigma^2}{|C_i|}
<\epsilon$, we can derive the necessary upper bound. Note that when $\theta$ is relatively small, this bound can be easily satisfied.

\end{proof}

\begin{proof}[Proof of Lemma \ref{lemma:1}]
\label{app:proof1}
	Since $\widehat A_t$ is estimated using data from $t+1$ to $t+w$, it is {\it independent} from $v_t$. This independence property allows for the straightforward computation of the two expectations in Lemma \ref{lemma:1} and contributes towards the proper selection of drift $d$. \textcolor{black}{Recall we denote $\hat{\Lambda}_{t} = \mathrm{diag}(\hat{\rho}_{t1},\ldots,\hat{\rho}_{tm})$ as the smallest $m$ eigenvalues of the sample covariance matrix, and $\hatU_t$ as the corresponding eigenvectors.}	Note that under the pre-change distribution we can write:
	\begin{align*}
		&\mathbb{E}_{\infty}[v_t^\top \widehat{A}_t\widehat{A}_t^\top v_t]\\=&\mathbb{E}_{\infty}[v_t^\top \hatU_t\hat{\Lambda}_{t}^{-1} \hatU_{t}^\top v_t ]
		=\mathbb{E}_{\infty}[v_t^\top (\sum_{i=1}^m\frac{1}{\hat{\rho}_{ti}} \hat{u}_{ti}\hat{u}_{ti}^\top) v_t ]\\
		=&\sum_{i=1}^m \mathbb{E}_{\infty}[\frac{1}{\hat{\rho}_{ti}} (\hat{u}_{ti}^\top v_t)^2]=\sum_{i=1}^m \mathbb{E}_{\infty}[\frac{1}{\hat{\rho}_{ti}} \hat{u}_{ti}^\top \mathbb{E}_{\infty}[v_t v_t^\top]\hat{u}_{ti}]\\
		=&\frac{1}{\sigma^2}\sum_{i=1}^m \mathbb{E}_{\infty}[\frac{1}{\hat{\rho}_{ti}} \hat{u}_{ti}^\top \hat{u}_{ti}]=\frac{1}{\sigma^2}\sum_{i=1}^m \mathbb{E}_{\infty}[\frac{1}{\hat{\rho}_{ti}} ] .
	\end{align*}
	
	Note that under the pre-change measure, the estimated $\hatU_t$ and $(\hatLambda_t)^{-1}$ are eigenvectors and eigenvalues for the sample covariance matrix constructed by samples generated from $\mathcal N(0,(1/\sigma^{2})I)$. 
	\textcolor{black}{We approximate the above expectation as $\mathbb{E}_{\infty}[1/\hat{\rho}_{ti} ]\approx \sigma^2$ using the eigenvalues of the ground truth sample covariance matrix. Note that such approximation is valid when the number of samples (i.e., the sliding window size $w$) is large enough, under such cases, the sample eigenvalues will highly concentrate around the true eigenvalues. As we will show later, the optimal window size $w$ is indeed sufficiently large (in the order of $\sqrt{\log\gamma}$) in the asymptotic regime we consider.}
	Following the discussions above, we write:
	\begin{equation}
		\mathbb{E}_{\infty}[v_t^\top \widehat{A}_t\widehat{A}_t^\top v_t] \approx
		m.
	\end{equation}
	
	

	Similarly, we can derive the results for post-change distribution:  
	\begin{align*}
		&\mathbb{E}_{0}[v_t^\top \widehat{A}_t\widehat{A}_t^\top v_t]\\
		=&\mathbb{E}_{0}[v_t^\top \hatU_t\Lambda_{t}^{-1} \hatU_{t}^\top v_t ]
		=\mathbb{E}_{0}[v_t^\top (\sum_{i=1}^m{\frac{1}{\hat\rho}_{ti}} \hat{u}_{ti}\hat{u}_{ti}^\top) v_t ]\\
		=&\sum_{i=1}^m \mathbb{E}_{0}[\frac{1}{\hat{\rho}_{ti}} (\hat{u}_{ti}^\top v_t)^2]
		=\sum_{i=1}^m \mathbb{E}_{0}[\frac{1}{\hat{\rho}_{ti}} \hat{u}_{ti}^\top \mathbb{E}_{0}[v_t v_t^\top]\hat{u}_{ti}]\\
 =& \textcolor{black}{\sum_{i=1}^m \mathbb{E}_{0}[\frac{1}{\hat{\rho}_{ti}} \hat{u}_{ti}^\top (\sum_{j=1}^n\rho_ju_ju_j^\top)\hat{u}_{ti}]}\\
		=&\sum_{i=1}^m\sum_{j=1}^n \mathbb{E}_{0}[\frac{\rho_j}{\hat{\rho}_{ti}} (\hat{u}_{ti}^\top {u}_{j})^2],
	\end{align*}
	where $1/\hat\rho_1,\ldots,1/\hat\rho_m$ are the estimated eigenvalues of $\widehat A\widehat A^\top$ and $\rho_1,\ldots,\rho_n$ are the eigenvalues of the true covariance matrix $(AA^T+\sigma^2 I)^{-1}$ \textcolor{black}{so that it can be written as $\mathbb{E}_{0}[v_tv_t^\top]=\sum_{j=1}^n\rho_ju_ju_j^\top$}.
As stated in Theorem \ref{thm1}, the estimated eigenvalue and eigenvector from the sample covariance matrix are independent, and Theorem \ref{thm1} characterizes the estimation error $ e_{it}=\hat{\varphi}_{ti}-u_i$, \textcolor{black}{which satisfies  $\sqrt{w}e_{it} \stackrel{d}{\rightarrow}\mathcal{N}\bigg(0, \sum_{k=1,k\neq i}^{m}\frac{\lambda_i \lambda_k}{(\lambda_i-\lambda_k)^2}u_ku_k^\top\bigg) $. Asymptotically, we have $\mathbb{E}_0[u_i^\top e_{it}]=u_i^\top\mathbb{E}_0[e_{it}]\rightarrow 0$ and
\begin{equation*}
    \mathbb{P}(u_i^\top e_{it}e_{it}^\top u_i > \epsilon)\leq \frac{\mathbb{E}_0[u_i^\top e_{it}e_{it}^\top u_i]}{\epsilon} =\frac{u_i^\top \text{Cov}(e_{it})u_i}{\epsilon}\rightarrow 0,
\end{equation*} 
thus $u_i$ is perpendicular to the error term $e_{it}$ with probability 1.}
	
Since the estimated eigenvector is normalized to be unit-norm, we have the following:
	\[
	\hat{u}_{it} =\frac{\hat{\varphi}_{ti}}{\left\Vert \hat{\varphi}_{ti}\right\Vert}= \frac{u_i + e_{it}}{\left\Vert u_i+e_{it}\right\Vert} .
	\]
	First of all, due to the orthogonality of different eigenvectors, we have:
	\[
	\sum_{i=1}^m\sum_{j=1}^n \mathbb{E}_{0}[\frac{\rho_j}{\hat{\rho}_{ti}} (\hat{u}_{ti}^\top {u}_{j})^2]
	=\sum_{i=1}^m \mathbb{E}_{0}[\frac{\rho_i}{\hat{\rho}_{ti}} (\hat{u}_{ti}^\top u_i)^2],
	\]
	where $u_i$, $i=1,\ldots,m$, are the eigenvectors corresponding to the smallest $m$ eigenvalues of the true covariance matrix.
	\textcolor{black}{Then we examine the term $\rho_i/\hat{\rho}_{ti}$, note that the true eigenvalue of the covariance matrix is $\rho_i=1/(|C_i|+\sigma^2) =1/(\lambda_i+\sigma^2)$ and the estimated eigenvalues of $(\widehat A\widehat A^\top)^{-1}$ is approximately $1/\lambda_i$. Thus we have $\mathbb{E}_{0}[\rho_i/\hat{\rho}_{ti}]=\lambda_i/(\lambda_i+\sigma^2)$, and of course such approximation is under the asymptotic case where the window size $w$ is sufficiently large.}
	
	Combining these together, we have
\textcolor{black}{
	\begin{align*}
		&\mathbb{E}_{0}[v_t^\top \widehat{A}_t\widehat{A}_t^\top v_t]\\
		=&\sum_{i=1}^m \mathbb{E}_{0}[\frac{\rho_i}{\hat{\rho}_{ti}} (\hat{u}_{ti}^\top u_i)^2]=
		\sum_{i=1}^m \frac{\lambda_i}{\lambda_i+\sigma^2}\mathbb{E}_0\bigg[\frac{(\left\Vert u_i \right\Vert^2+e_{it}^\top u_i)^2}{\left\Vert u_i+e_{it} \right\Vert^2}\bigg]\\
		=&  \sum_{i=1}^m \frac{\lambda_i}{\lambda_i+\sigma^2} \mathbb{E}_0\bigg[\frac{(1+e_{it}^\top u_i)^2}{1+\left\Vert e_{it} \right\Vert^2}\bigg]
		=  \sum_{i=1}^m \frac{\lambda_i}{\lambda_i+\sigma^2} \mathbb{E}_0\bigg[\frac{1}{1+\left\Vert e_{it}\right\Vert^2}\bigg]\\
  =&\textcolor{black}{
    \sum_{i=1}^m \frac{\lambda_i}{\lambda_i+\sigma^2} \mathbb{E}_0 \bigg[\frac{1-\Vert e_{it} \Vert^4+\Vert e_{it} \Vert^4}{1+\left\Vert e_{it} \right\Vert^2}\bigg] 
   } \\
   =&\textcolor{black}{
    \sum_{i=1}^m \frac{\lambda_i}{\lambda_i+\sigma^2} \mathbb{E}_0 \bigg[\frac{(1+\Vert e_{it} \Vert^2)(1-\Vert e_{it} \Vert^2)+\Vert e_{it} \Vert^4}{1+\left\Vert e_{it} \right\Vert^2}\bigg] 
   } \\
	=&  \sum_{i=1}^m \frac{\lambda_i}{\lambda_i+\sigma^2} \bigg\{1-\mathbb{E}_0[\left\Vert e_{it} \right\Vert^2]+\mathbb{E}_0\bigg[\frac{\left\Vert e_{it}\right\Vert^4}{1+\left\Vert e_{it}\right\Vert^2}\bigg]\bigg\}.
	\end{align*}}
 \textcolor{black}{The second equality is by substituting $\hat{u}_{it} = (u_i + e_{it})/{\left\Vert u_i+e_{it}\right\Vert}$, and the third equality is due to the fact that $u_i$ is an unit-norm vector and $u_i$ is perpendicular to $e_{it}$ with probability 1. 
}
	For the two expectations above, using Gaussian approximation from Theorem \ref{thm1}, we have
	\[
	\mathbb{E}_0[\left\Vert e_{it} \right\Vert^2] = \frac{1}{w} \sum_{k=1,k\neq i}^{m}\frac{\lambda_i \lambda_k}{(\lambda_i-\lambda_k)^2}=\frac{B_i}{w},
	\]
	and \textcolor{black}{to estimate $\mathbb{E}_0\big[{\left\Vert e_{it}\right\Vert^4}/{(1+\left\Vert e_{it}\right\Vert^2)}\big]$, we provide the following lemma:
	\begin{lemma}
	\label{lemma:forth-moment}
	For a multivariate Gaussian random variable $x\sim \mathcal{N}(0,\Sigma)$ where $\Sigma_{ij}=\sigma_{ij}^2$, we have
	    $$\mathbb{E}(\Vert x\Vert_2^4)\leq 3(\mathrm{tr}(\Sigma))^2.$$
     	\end{lemma}
	We can derive this lemma using the properties of uni-variate Gaussian distribution where $x_i\sim \mathcal{N}(0,\sigma_{ii}^2)$ $\mathbb{E}(x_i^4)=3\sigma_{ii}^4$  and Cauchy-Schwarz inequality (details omitted).
	As a result, according to Lemma \ref{lemma:forth-moment}, we can conclude that random variable $e_{it}$ that follows a Gaussian distribution have $\mathbb{E}_0[\left\Vert e_{it}\right\Vert^4]=\mathbb{E}_0[(\left\Vert e_{it}\right\Vert^2)^2] \leq 3 [\mathrm{tr}(\text{Cov}(e_{it}))]^2={3B_i^2}/{w^2}$. Then we have:}
	\begin{align*}
		0&\leq \mathbb{E}_0\bigg[\frac{\left\Vert e_{it}\right\Vert^4}{1+\left\Vert e_{it}\right\Vert^2}\bigg]\\
		 &\leq \mathbb{E}_0[\left\Vert e_{it}\right\Vert^4] \leq  \frac{3}{w^2} \bigg(\sum_{k=1,k\neq i}^{m}\frac{\lambda_i \lambda_k}{(\lambda_i-\lambda_k)^2}\bigg)^2
		=\frac{3B_i^2}{w^2}.
	\end{align*}

	Therefore, we get the desired approximation to the post-change expectation: 
	\begin{align*}
	\widetilde D&=\sum_{i=1}^m \frac{\lambda_i}{\lambda_i+\sigma^2} \bigg(1-\frac{B_i}{w}\bigg) \leq \mathbb{E}_{0}[v_t^\top \widehat{A}_t\widehat{A}_t^\top v_t]\\
	&\leq  \sum_{i=1}^m \frac{\lambda_i}{\lambda_i+\sigma^2} \bigg(1-\frac{B_i}{w}+\frac{3B_i^2}{w^2}\bigg) = D.
	\end{align*}
\textcolor{black}{From above, we see that when the window size $w$ becomes large, $\widetilde D/D\rightarrow 1$, we then have $\mathbb{E}_{0}[v_t^\top \widehat{A}_t\widehat{A}_t^\top v_t]$ can be well approximated by $D$.}	
	
\end{proof}

\begin{proof}[Proof of Lemma \ref{lemma:KL}]

By definition of KL divergence for emerging subspace case, \begin{align*}
\mathcal{I}_0&=\mathbb{E}_0\bigg[\log\frac{f_0(v_t)}{f_\infty(v_t)}\bigg]\\
&=\mathbb{E}_0\bigg[-\frac{1}{2}v_t^\top AA^\top v_t+\frac{1}{2}\log\frac{\mbox{det}(AA^\top+\sigma^2 I)}{\sigma^{2n}}\bigg]\\
&=-\frac{1}{2}\mathbb{E}_0\bigg[v_t^\top AA^\top v_t\bigg]+\frac{1}{2}\log\frac{\mbox{det}(AA^\top+\sigma^2 I)}{\sigma^{2n}}.
\end{align*}
Since A is not a full-rank matrix, the eigendecomposition of $AA^\top=U\Lambda U^\top$ has dimension $U\in\mathbb{R}^{n\times m}$ and $\Lambda \in \mathbb{R}^{m\times m}$. By adding basis $(\tilde{U}\in\mathbb{R}^{(n\times(n-m))})$ from nullspace to make $U'$ to be a square matrix, we have the following:
\[
U' =
\begin{bmatrix}
	U & \tilde{U} \\
\end{bmatrix}\in\mathbb{R}^{n\times n},
\Lambda' =
\begin{bmatrix}
	\lambda_1 & &&& \\
	& \ddots &&& \\
	& & \lambda_m&&\\
	&&& 0&&\\
	&&&&\ddots&\\
	&&&&&0
\end{bmatrix}\in\mathbb{R}^{n\times n},
\]
where $U'$ is an orthogonal matrix of rank $n$ where $U'^\top U'=I$ and  $AA^\top=U'\Lambda'U'^\top$. Then the first term of $\mathcal{I}_0$ can be derived as:
\begin{align*}
&\mathbb{E}_0\bigg[v_t^\top AA^\top v_t\bigg]\\=&\mathbb{E}_0\bigg[\text{tr}(AA^\top v_t v_t^\top )\bigg]
=\text{tr}\bigg[AA^\top\mathbb{E}_0(vv^\top)\bigg]\\
=&\text{tr}\bigg[AA^\top (AA^\top+\sigma^2 I)^{-1} \bigg]=\text{tr}\bigg[U'\Lambda' U'^\top (U'(\Lambda'+\sigma^2 I)^{-1}U'^\top) \bigg]\\
=&\text{tr}\bigg[U'\Lambda'(\Lambda'+\sigma^2I)^{-1}U'^\top\bigg]=\text{tr}\bigg[\Lambda'(\Lambda'+\sigma^2I)^{-1}U'^\top U'\bigg]\\
=&\text{tr}\bigg[\Lambda'(\Lambda'+\sigma^2I)^{-1}\bigg]=\sum_{i=1}^m \frac{\lambda_i}{\sigma^2+\lambda_i}.
\end{align*}
Moreover, the second term becomes
\begin{align*}
\frac{1}{2}\log\frac{\mbox{det}(AA^\top+\sigma^2 I)}{\sigma^{2n}}
&= \frac 1 2 \sum_{i=1}^m \log
\left(
\frac{\lambda_i}{\sigma^2}  +1
\right) \\
&= -\frac 1 2 \sum_{i=1}^m \log
\left( 1-
\frac{\lambda_i}{\sigma^2+\lambda_i} 
\right). 
\end{align*}
Thus we have shown that:
\[
    \mathcal{I}_0 =
    -\frac 1 2 \sum_{i=1}^m h\left(
    \frac{\lambda_i}{\sigma^2+\lambda_i} 
    \right),
    \]
    where $h(x) = x + \log(1-x)$. 
\end{proof}

\begin{proof}[Proof of Lemma \ref{lemma:2}]
Following previous results, we define $\psi(w)= \mathbb{E}_0(\mathcal{T}_C)$ in \eqref{eq:EDD_min} as a function of $w$ and try to find the $w^*$ that minimizes $\psi(w)$. Recall we have
\begin{align*}
\psi(w) &= \frac{2\log\gamma\big(1+o(1)\big)}{-2\delta_\infty \mathbb{E}_0[v_t^\top \hatU_t\hatLambda^{-1}_t\hatU_t^\top v_t]  + m\log ( 1+2\delta_\infty)}+w \\
&= \frac{2\log\gamma\big(1+o(1)\big)}{-2\delta_\infty D  + m\log ( 1+2\delta_\infty)}+w.
\end{align*}
%
Substitute $D$ into $\psi(w)$ gives us the expression for the EDD as:
	\[
	\psi(w) =\frac{2\log\gamma}{D-m+m\log(m/D)}+w = \frac{2\log\gamma}{mg(D/m)} + w.
	\]
	Recall we have denoted $\Delta =\frac{1}{m}\sum_{i=1}^m\frac{\lambda_i}{\sigma^2+\lambda_i}$, thus we have
	\begin{align*}
	D&= \sum_{i=1}^m\frac{\lambda_i}{\sigma^2+\lambda_i} -\sum_{i=1}^m\frac{\lambda_iB_i}{\sigma^2+\lambda_i}\frac{1}{w}+\sum_{i=1}^m\frac{3B_i^2\lambda_i}{\sigma^2+\lambda^2}\frac{1}{w^2}\\&=m\Delta-\sum_{i=1}^m\frac{\lambda_iB_i}{\sigma^2+\lambda_i}\frac{1}{w}+\sum_{i=1}^m\frac{3B_i^2\lambda_i}{\sigma^2+\lambda^2}\frac{1}{w^2}.
	\end{align*}
Note that $m\Delta \in (0,m)$ and it is close to close to 0 when $\sigma^2$ is large and close to $m$ when $\sigma^2$ is small. 
Let 
	\begin{align*}
		B_i'&=\frac{\lambda_iB_i}{\sigma^2+\lambda_i};\quad C_i'=\frac{3B_i^2\lambda_i}{\sigma^2+\lambda^2}.
	\end{align*}
	We have 
	\[
	D=m\Delta-\frac{\sum B_i'}{w}+\frac{\sum C_i'}{w^2}.
	\]
	Then
	\begin{align*}
		&\psi(w) = w+\\
		&\frac{2\log\gamma}{m(\Delta-1)-\frac{\sum B_i'}{w} + \frac{\sum C_i'}{w^2} -m\log(\Delta-\frac{\sum B_i'}{mw}+\frac{\sum C_i'}{mw^2})}.
	\end{align*}
	Note that 
	\begin{align*}
		&\Delta - \frac{\sum B_i'}{mw}+\frac{\sum C_i'}{mw^2}=\Delta(1- \frac{\sum B_i'}{\Delta mw}+\frac{\sum C_i'}{\Delta mw^2}).
	\end{align*}
	Apply the Taylor expansion for $\log(1-x)$, we have
	\begin{align*}
		&\log\bigg(1- \frac{\sum B_i'}{\Delta mw}+\frac{\sum C_i'}{\Delta mw^2}\bigg)=  - \frac{\sum B_i'}{\Delta mw}+\frac{\sum C_i'}{\Delta mw^2} + o(1/w).
	\end{align*}
	Substitute into the EDD $\psi(w)$ we can obtain that:
	\begin{align*}
	&\psi(w)=w+
	(2\log \gamma)\bigg/\bigg[m(\Delta-1)-m\log(\Delta)\\
	&-\frac{1}{w}(1-\frac{1}{\Delta})(\sum_i B_i')+\frac{1}{w^2}(1-\frac{1}{\Delta}))(\sum_i C_i')\bigg],
	\end{align*}
	and note that the denominator will be dominated by the term $m(\Delta-1)-m\log(\Delta)$. Let the first-order derivative of $\psi(w)$ equal to 0 and ignore $o(1/w)$ terms, we obtain
	\begin{align*}
	\psi'(w)&=1+\frac{2(\log\gamma)(\frac{1}{\Delta}-1) (\sum B_i')}{[m(\Delta-1)-m\log(\Delta)-(\frac{1}{\Delta}-1) (\sum B_i')/w]^2 \cdot w^2}\\
	&\approx 1-\frac{2(\log\gamma)(\frac{1}{\Delta}-1) (\sum B_i')}{[m(\Delta-1)-m\log(\Delta)]^2 \cdot w^2}. 
	\end{align*}
	The derivation above is approximated by ignoring $o(1/w)$ term in the denominator in the second step. It can be seen clearly that when $\psi '(w^*)=0$ we get the minimum value for $\psi(w)$. Thus we have:
	\begin{align*}
		w^* &= \frac{\sqrt{2(\log \gamma) (\frac{1}{\Delta}-1) (\sum B_i')}}{m(\Delta-1)-m\log(\Delta)}\\
		&= \frac{\sqrt{2(\log \gamma) (\frac{1}{\Delta}-1)(\sum \frac{\lambda_i B_i}{\sigma^2+\lambda_i})}}{m(\Delta-1-\log(\Delta))}\\
		&=\frac{\sqrt{2(\log \gamma) (\frac{1}{\Delta}-1)(\sum_{i=1}^m \frac{\lambda_i B_i}{\sigma^2+\lambda_i})}}{mg\left(\Delta\right)}.
	\end{align*}
	Recall $g(x) =x-1-\log(x)$.
	Note that the denominator is always positive since $g(x)>0$ for all $x\in(0,1)$ and $\Delta\in(0,1)$.

\label{app:proof2}
\end{proof}

\begin{proof}[Proof of Theorem \ref{thm:main}]
\label{app:proof_thm}
\textcolor{black}{
To explain the $o(1)$ terms in \eqref{eq:thmcusum}, note that the result in \cite{siegmund2013sequential} states $\mathbb{E}_\infty[T_C]=e^b[1 - (b + 1)e^{-b}]/\mathcal{I}_\infty$, and $\mathbb{E}_0 [T_C]= b[1+(e^{-b}-1)b^{-1}]/\mathcal{I}_0$, and thus the two $o(1)$ terms are on the order of $be^{-b}$ and $b^{-1}$ respectively, or equivalently, the two $o(1)$ terms are on the order of $(\log\gamma)/\gamma$ and $1/\log\gamma$ respectively.}
Now using the more precise version of the EDD for the exact CUSUM \eqref{exact_CUSUM} and the Spectral CUSUM \eqref{eq:EDD_min}
\begin{align*}
\mathbb{E}_0[T_C]&=\frac{b (1+\mathcal O(b^{-1}))}{\mathcal I_0}. \quad \mbox{(exact CUSUM)}\\
\mathbb{E}_0(\mathcal{T}_C)&=\frac{2b\big(1+ \mathcal O(b^{-1})\big)}{m g(D/m)}+w. \quad \mbox{(Spectral CUSUM)}
\end{align*}
Recall $g(x) =x-1-\log(x)$.

Since, to achieve ARL constraint $\gamma$, the threshold $b$ for both procedures are on the order of $\log \gamma$, we obtain the ratio:
\begin{align*}
    \frac{\mathbb{E}_0[\mathcal{T}_C]}{\mathbb{E}_0[T_C]} &=
    \frac{2\mathcal I_0}{m g(D/m)}+
    \frac{w \mathcal I_0}{b(1+\mathcal O(b^{-1})}
    \\
    &=
     \frac{2\mathcal I_0/m}{ g(D/m)}+
    \frac{w \mathcal I_0}{b}(1+\mathcal O(b^{-1})),
\end{align*} 
using Taylor expansion $1/(1+x) = 1-x+x^2/2 \cdots$.
%
%
Note that
\[
    \frac{2\mathcal{I}_0}{m} =
    -\frac{1}{m}\sum_{i=1}^m h\left(
    \frac{\lambda_i}{\sigma^2+\lambda_i} 
    \right)= - \Delta -\frac {1}{m}\sum_{i=1}^m \log\left(
    \frac{\sigma^2}{\sigma^2+\lambda_i}\right), \] where we recall that $h(x) = x+\log(1-x)$,  $\Delta=\frac1m\sum_{i=1}^m 
    \frac{\lambda_i}{\sigma^2+\lambda_i}$. Meanwhile we can write
    \begin{align*}
    D&= \sum_{i=1}^m \frac{\lambda_i}{\lambda_i+\sigma^2} \bigg(1-\frac{B_i}{w}+\frac{3B_i^2}{w^2}\bigg) \\
    &= m\Delta -\sum_{i=1}^m\frac{\lambda_iB_i}{\sigma^2+\lambda_i}\frac{1}{w}+\sum_{i=1}^m\frac{3B_i^2\lambda_i}{\sigma^2+\lambda^2}\frac{1}{w^2},
    \end{align*} and thus 
    \[g(D/m) = g\left(\Delta -\sum_{i=1}^m\frac{\lambda_iB_i}{\sigma^2+\lambda_i}\frac{1}{mw}+\sum_{i=1}^m\frac{3B_i^2\lambda_i}{\sigma^2+\lambda^2}\frac{1}{mw^2} \right).\] When $w$ is sufficiently large (e.g., in our setting $w=\sqrt{\log\gamma}$, $\gamma \rightarrow\infty$), due to monotonicity and continuity of $g$, we have $g(D/m) \rightarrow g(\Delta) = \Delta - 1 - \log(\Delta)$ and 
    \begin{equation}
     \frac{2\mathcal I_0/m}{ g(D/m)}
    \rightarrow  \frac{\Delta +\frac {1}{m}\sum_{i=1}^m \log\left(
    \frac{\sigma^2}{\sigma^2+\lambda_i}\right) }{1-\Delta +  \log(\Delta)} = \mathcal O(1),
    \label{eq:ratio}
    \end{equation}
    (which does not depend on $w$ or $\gamma$) since $g(\Delta)>0$ and $\Delta \in (0, 1)$ for $\sigma^2$ bounded away from 0.
    Combining with $b$ is on the order of $\log \gamma$, $w$ is in the order of $\sqrt{\log \gamma}$, we are done.
    
    We can further evaluate the limiting upper bound, using the concavity of the logarithm function, \eqref{eq:ratio} can be further bounded by 
    \[
    \frac{\Delta +\frac {1}{m}\sum_{i=1}^m \log\left(
    \frac{\sigma^2}{\sigma^2+\lambda_i}\right) }{1-\Delta +  \log(\Delta)}
    \leq 
    \frac{\Delta + \log(1-\Delta)}{1-\Delta + \log \Delta},
    \]
    and the figure below plots the upper bound for $\Delta \in (0.4, 0.9)$.
    
    \begin{figure}
        \centering
        \includegraphics[width = 0.5\textwidth]{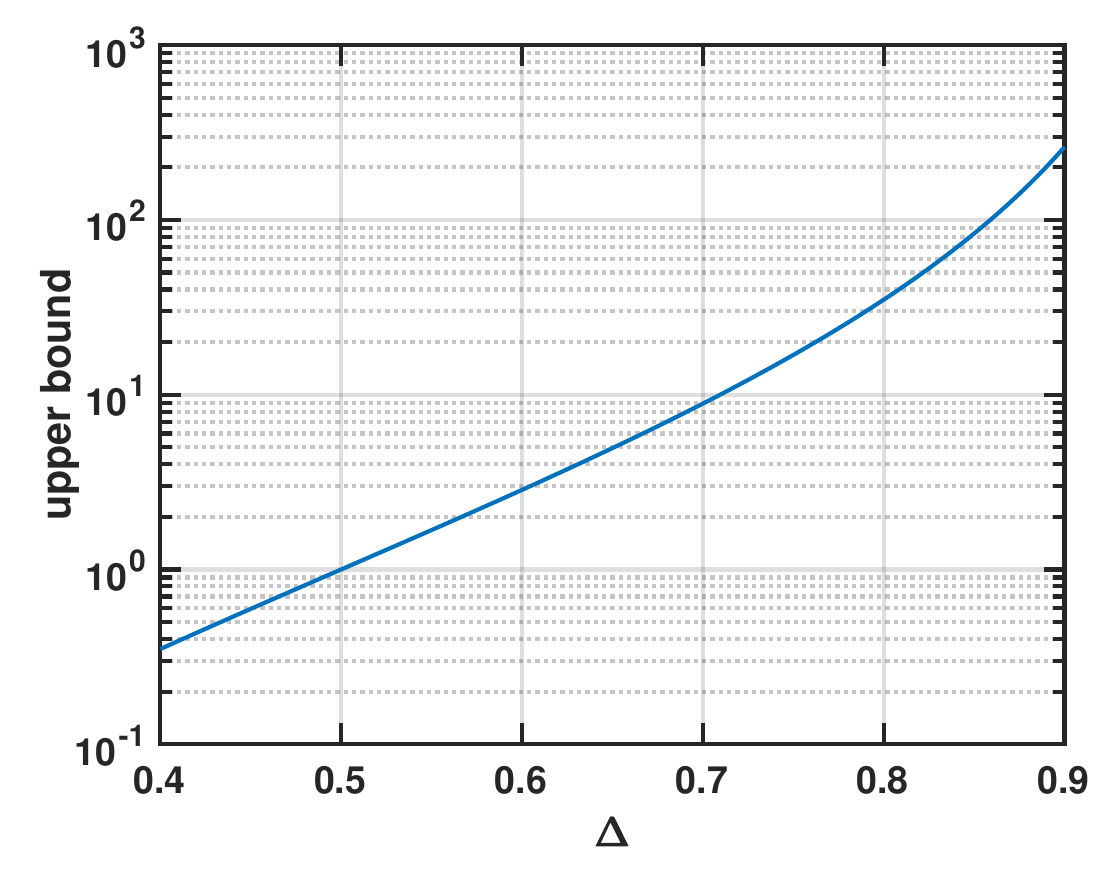}
        \caption{Plot of the asymptotic upper bound $\frac{\Delta + \log(1-\Delta)}{1-\Delta + \log \Delta}$ as a function $\Delta$, for $\Delta \in (0.4, 0.9)$. Note that it is reasonable for $\Delta$ bounded away from 1 (i.e., when $\sigma^2$ is bounded away from 0.)}
        \label{fig:my_label}
    \end{figure}
    
\end{proof}
\end{document}